\documentclass[12pt]{article}
\usepackage{amsfonts}
\usepackage{amsmath}
\usepackage{amssymb}
\usepackage{amsthm}
\usepackage[numbers]{natbib}
\usepackage{booktabs}
\usepackage{algorithm}
\usepackage{algorithmic}
\usepackage{pgfplots}
\pgfplotsset{compat=1.8}
\usepackage{graphicx}
\usepackage{color}
\usepackage{colortbl}   
\usepackage{a4wide}
\usepackage{multirow}
\usepackage{boxedminipage}
\usepackage{enumerate}
\usepackage{subfig}
\usepackage{textcomp}
\usepackage{hyperref}
\usepackage{placeins}

\newcommand{\N}{\ensuremath{\mathbb{N}}}

\newcommand{\T}{\ensuremath{\mathbb{T}}}

\newcommand{\Z}{\ensuremath{\mathbb{Z}}}

\newcommand{\R}{\ensuremath{\mathbb{R}}}
\newcommand{\C}{\ensuremath{\mathbb{C}}}

\newcommand{\ii}{\textnormal{i}}
\newcommand{\e}{\textnormal{e}}

\newcommand{\ceil}[1]{\left\lceil#1\right\rceil}

\newcommand{\zb}[1]{\ensuremath{\boldsymbol{#1}}}

\renewcommand{\ln}{\mathrm{ln\,}}
\newcommand{\bigtimes}{\mathop{\text{\Large{$\times$}}}}
\newcommand{\boldx}{{\ensuremath{\boldsymbol{x}}}}
\newcommand{\boldk}{{\ensuremath{\boldsymbol{k}}}}
\newcommand{\boldh}{{\ensuremath{\boldsymbol{h}}}}
\newcommand{\boldl}{{\ensuremath{\boldsymbol{l}}}}
\newcommand{\boldz}{{\ensuremath{\boldsymbol{z}}}}

\newcommand{\boldj}{{\ensuremath{\boldsymbol{j}}}}

\newcommand{\boldzero}{{\ensuremath{\boldsymbol{0}}}}
\newcommand{\boldone}{{\ensuremath{\boldsymbol{1}}}}

\definecolor{darkgreen}{rgb}{0.0,0.5,0.0}
\definecolor{darkorange}{RGB}{255,90,0}

\DeclareMathOperator{\ml}{mr1l}
\DeclareMathOperator{\sgn}{sgn}

\newtheorem{theorem}{Theorem}[section]
\newtheorem{lemma}[theorem]{Lemma}
\newtheorem{remark}[theorem]{Remark}
\newtheorem{generalisation}[theorem]{Generalisation}
\newtheorem{definition}[theorem]{Definition}
\newtheorem{example}[theorem]{Example}
\newtheorem{corollary}[theorem]{Corollary}
\newtheorem{proposition}[theorem]{Proposition}

\newenvironment{Theorem}{\goodbreak \begin{theorem}\sl}{\end{theorem}}
\newenvironment{Lemma}{\goodbreak \begin{lemma}\sl}{\end{lemma}}
\newenvironment{Remark}{\goodbreak \begin{remark}\rm}{\bend\end{remark}}

\newenvironment{Example}{\goodbreak \begin{example}\rm}{\bend\end{example}}

\newenvironment{Corollary}{\goodbreak \begin{corollary}\sl}{\end{corollary}}

\makeatletter
\def\imod#1{\allowbreak\mkern10mu({\operator@font mod}\,\,#1)}
\makeatother

\pgfdeclareplotmark{hexagon}
{%
  \pgfpathmoveto{\pgfqpoint{0pt}{\pgfplotmarksize}}
  \pgfpathlineto{\pgfqpointpolar{150}{\pgfplotmarksize}}
  \pgfpathlineto{\pgfqpointpolar{210}{\pgfplotmarksize}}
  \pgfpathlineto{\pgfqpointpolar{270}{\pgfplotmarksize}}
  \pgfpathlineto{\pgfqpointpolar{330}{\pgfplotmarksize}}
  \pgfpathlineto{\pgfqpointpolar{30}{\pgfplotmarksize}}
  \pgfpathclose
\pgfusepathqstroke
}

\pgfdeclareplotmark{hexagon*}
{%
  \pgfpathmoveto{\pgfqpoint{0pt}{\pgfplotmarksize}}
  \pgfpathlineto{\pgfqpointpolar{150}{\pgfplotmarksize}}
  \pgfpathlineto{\pgfqpointpolar{210}{\pgfplotmarksize}}
  \pgfpathlineto{\pgfqpointpolar{270}{\pgfplotmarksize}}
  \pgfpathlineto{\pgfqpointpolar{330}{\pgfplotmarksize}}
  \pgfpathlineto{\pgfqpointpolar{30}{\pgfplotmarksize}}
  \pgfpathclose
\pgfusepathqfillstroke
}

\numberwithin{equation}{section}
\numberwithin{table}{section}
\numberwithin{figure}{section}

\newcommand{\bend}{\hspace*{0ex} \hfill \hbox{\vrule height
    1.5ex\vbox{\hrule width 1.4ex \vskip 1.4ex\hrule  width 1.4ex}\vrule
    height 1.5ex}}

\long\def\symbolfootnote[#1]#2{\begingroup%
\def\thefootnote{\fnsymbol{footnote}}\footnote[#1]{#2}\endgroup}
\newcommand{\sspan}{\textnormal{span}}

\newcommand{\OO}[1]{\mathcal{O}\left(#1\right)}

\renewcommand{\mathbf}[1]{\ensuremath{\boldsymbol{#1}}}
\renewcommand{\textbf}[1]{{\ensuremath{\boldsymbol{#1}}}}
\renewcommand{\thefootnote}{\fnsymbol{footnote}}

\title{Approximation of multivariate periodic functions based on sampling along multiple rank-1 lattices}
\date{}
\author{Lutz K\"ammerer\footnotemark[1]
        \and Toni Volkmer\footnotemark[2]}
        
\hypersetup{
          pdfauthor=Lutz K\"{a}mmerer and Toni Volkmer,
          pdftitle={Approximation of multivariate periodic functions based on sampling along multiple rank-1 lattices},
          pdfsubject={65T40, 42B05, 68Q25, 42B35, 65T50, 65Y20, 68Q17, 65D30, 65D32},
          pdfkeywords={approximation of multivariate periodic functions, trigonometric polynomials, generalized hyperbolic cross, generalized mixed smoothness, lattice rule, multiple rank-1 lattice, fast Fourier transform},
          pdfcreator = {pdflatex},
          plainpages=false,
          pdfstartview=Fit,       
          pdfpagemode=UseOutlines, 
          bookmarksnumbered=true, 
          bookmarksopen=false,     
          bookmarksopenlevel=0,   
          colorlinks=true,       
          linkcolor=black,
          citecolor=black,
          urlcolor=black}

\begin{document}

\maketitle

\footnotetext[1]{\scriptsize
  Chemnitz University of Technology, Faculty of Mathematics, 09107 Chemnitz, Germany\\
  lutz.kaemmerer@mathematik.tu-chemnitz.de, Phone:+49-371-531-37728, Fax:+49-371-531-837728
}
\footnotetext[2]{\scriptsize
  Chemnitz University of Technology, Faculty of Mathematics, 09107 Chemnitz, Germany\\
  toni.volkmer@mathematik.tu-chemnitz.de, Phone:+49-371-531-39999, Fax:+49-371-531-839999
}

\begin{abstract}
In this work, we consider the approximate reconstruction of high-dimensional periodic functions based on sampling values.
As sampling schemes, we utilize so-called reconstructing multiple rank-1 lattices, which combine several preferable properties such as
easy constructability, the existence of high-dimensional fast Fourier transform algorithms, high reliability, and low oversampling factors.
Especially, we show error estimates
for functions from Sobolev Hilbert spaces of generalized mixed smoothness.
For instance, when measuring the sampling error in the $L_2$-norm,
we show sampling error estimates where the exponent of the main part reaches those of the optimal sampling rate except for an offset of $1/2+\varepsilon$,
i.e., the exponent is almost a factor of two better up to the mentioned offset compared to single rank-1 lattice sampling.
Various numerical tests in medium and high dimensions demonstrate the high performance and confirm the obtained theoretical results of multiple rank-1 lattice sampling.
\end{abstract}



\section{Introduction}

High-dimensional integrals are often treated numerically
by applying a cubature formula. 
Commonly used methods are equally weighted cubature formulas,
called quasi-Monte Carlo rules, which average
a specific set of sampling values of the integrand.
One established type of these methods are
lattice rules, where the sampling nodes have 
a group structure, cf.\ \cite{SlJo94, JoKuSl13}.
In particular, the simple structure of so-called rank-1 lattices allows for theoretical analysis,
which was already exploited in early contributions, see e.g.\ \cite{Ko59, Zaremba72}.
Later, component--by--component construction approaches rediscovered in \cite{SlRe02}
distinctly improved the applicability and shifted rank-1 lattice rules
into the focus of research, again.

Furthermore, these lattice rules have been widely used
for the approximation of high-dimensional periodic functions by
approximately computing a subset of the Fourier coefficients
\begin{equation}
\hat{f}_{\zb k}:=\int_{\T^d}f(\zb x) \; \e^{-2\pi\ii\zb k\cdot\zb x}\mathrm{d}\zb x, \quad \boldk\in\Z^d,\label{equ:fc:integral}
\end{equation}
 of a suitable periodic function $f\colon \T^d \rightarrow \C$ based on sampling values, where $\T^d\simeq [0,1)^d$ is the $d$-dimensional torus.
Due to the preferable properties of rank-1 lattices mentioned above, many contributions 
estimate approximation errors for rank-1 lattice sampling,
cf.\ e.g.\ \cite{Ko63, Tem86, Tem93, KuSlWo08, KuWaWo09, KaPoVo13, KaPoVo14, ByKaUlVo16}.
We stress on the fact that these theoretical results may be 
of particular relevance in practice since the structure of lattices allows for a very efficient
simultaneous computation of Fourier coefficients using fast Fourier transforms (FFTs), cf.\ \cite{LiHi03, Be13}.
The arithmetic complexity of these FFTs is almost linear in the amount of data, i.e. the number of sampling values, that has to be handled.
This is one of the reasons why recent research focuses on estimates on the relation of the number of used sampling values
to the approximation errors.

One of the most recent results for rank-1 lattices is a negative one, i.e., 
the relation of the number of used sampling values
to the approximation error is far away from the optimum
in highly interesting approximation settings,  cf.\ \cite{ByKaUlVo16}.
Moreover, approximation errors in a best possible order of magnitude
require specific rank-1 lattices that need to be determined.
Different component--by--component methods allow for
the construction of suitable rank-1 lattices, cf.\ e.g.\ \cite{KuSlWo06, CoNu06, Kae2012}.
However, these construction methods require
a number of arithmetic operations that is at least linear in the
number of sampling nodes within the rank-1 lattice, and thus,
also suffer from the necessarily huge number of sampling nodes.

In order to avoid this huge oversampling, which is caused by the structure of single rank-1 lattices,
a new type of spatial discretizations for multivariate trigonometric polynomials
was developed recently, cf.\ \cite{Kae16}.
Roughly speaking, one joins multiple rank-1 lattices in order to determine
a sampling scheme. Exploiting the structure of each of the joined rank-1 lattices,
a fast Fourier transform algorithm arises in a natural way.

Based on this specific type of spatial discretizations for
multivariate tri\-go\-nometric polynomials, we investigate the 
arising sampling method as a sampling operator for functions belonging to specific function spaces.
In this paper, we specify first approximation properties of this new sampling method,
that is based on the approximation of Fourier coefficients using a set of
rank-1 lattice rules that fulfills certain properties.
The construction of those sampling sets is easy
and very efficient, cf.~\cite{Kae17}.
Due to the structure of the sampling sets,
the computation of the aforementioned approximants is eminently efficient and
also simple, cf.\ Algorithm~\ref{algo:mlfft2:ifft_direct}, originally stated as \cite[Algorithm~2]{KaPoVo17}.
Moreover, numerical tests indicate that the Fourier transform is stable, cf.~\cite{Kae16,Kae17}.

The main focus of our considerations is on the decay
of worst case approximation errors for sampling methods (worst case sampling errors)
for increasing numbers of function samples.
For a sampling set $\mathcal{G} := \{\boldx_1,\ldots,\boldx_M\}$, $M\in\N$, normed function spaces $\mathcal{F}$ and $Y$, we denote
the best possible worst case sampling errors measured in the norm of the target space $Y$
for functions belonging to the source space $\mathcal{F}$
as
\begin{align*}
\operatorname{Samp}_{\mathcal{G}}({\mathcal{F}},Y):=\inf_{A:\C^M\to Y}\sup_{\|f|{\mathcal{F}}\|\leq 1}\left\|f-A\big(f(\boldx_j)\big)_{j=1}^M\right\|_Y\,.
\end{align*}
Here, $A$ denotes possibly non-linear sampling operators that
constitute an approximation of $f\in\mathcal{F}$ in $Y$ using at most $M$
sampling values $f(\boldx_j), \boldx_j\in\mathcal{G}$.
Then, we define the general sampling numbers by
$$
g_M({\mathcal{F}},Y):=\inf_{|\mathcal{G}|\le M}\operatorname{Samp}_
{\mathcal{G}}({\mathcal{F}},Y)\,,\quad M\in \N,
$$
which is the best possible worst case sampling error one can achieve using at most
$M$ sampling values. Additional restrictions on the algorithms or the sampling sets will
result in specific sampling numbers, that are at least as big as $g_M({\mathcal{F}},Y)$.

\sloppy
In particular, the sampling number for linear operators, that we will denote by $g_M^{\mathrm{lin}}({\mathcal{F}},Y)$,
can be determined by restricting the operators $A$ to linear ones.
For specific choices of source and target spaces the sampling numbers $g_M^{\mathrm{lin}}({\mathcal{F}},Y)$
are well known up to some logarithmic gaps, cf.\ column four in Tables \ref{tab:mainres:overview_general_r1l} and \ref{tab:mainres:overview_general_r1l_2},
and hence, are used for comparison.
For an overview on this topic we refer to \cite{ByDuSiUl14, DuTeUl16} and the references therein.

In this paper, we deal with a fixed structure of the sampling sets
\begin{equation}
\mathcal{G}=\Lambda:=\Lambda(\boldz_1,M_1,\ldots,\boldz_L,M_L):=
\bigcup_{\ell=1}^L \Lambda(\boldz_\ell,M_\ell),\label{eq:def:mr1l}
\end{equation}
so-called multiple rank-1 lattices, cf.~\cite{Kae16,Kae17}, which are the union of
single rank-1 lattices
\begin{equation}
\Lambda(\boldz_\ell,M_\ell):=\Big\{\frac {j}{M_\ell}\boldz_\ell\bmod \boldone\colon j=0,\ldots,M_\ell-1 \Big\} \subset \T^d, \quad \boldz_\ell\in\Z^d, \; M_\ell\in\N.
\label{eq:def:r1l}
\end{equation}
To this end, we define the corresponding sampling numbers
\begin{align*}
g_M^{\ml} ({\mathcal{F}},Y):=\inf_{|\Lambda|\le M}\operatorname{Samp}_
{\Lambda}({\mathcal{F}},Y)\,,\quad M\in \N,
\end{align*}
for sampling sets consisting of multiple rank-1 lattices~$\Lambda$.
Furthermore, we analyze a fixed linear algorithm $A^{\ml}$, cf.\ Algorithm~\ref{algo:mlfft2:ifft_direct}.
For that reason, the results of this paper are upper bounds on the number
$g_M^{\ml} ({\mathcal{F}},Y)$.
An additional restriction on the sampling sets $\Lambda$, i.e., $L=1$, leads to the definition
of the sampling numbers
\begin{align*}
g^{\mathrm{latt}_1}_M({\mathcal{F}},Y):=\inf_{\boldz\in\Z^d}\operatorname{Samp}_{\Lambda(\boldz,M)}({\mathcal{F}},Y)\,,\quad M\in \N,
\end{align*}
for sampling sets $\mathcal{G}$ that are exactly one rank-1 lattice. These
sampling numbers are already investigated in \cite{ByKaUlVo16} and 
serve for comparison to the results for multiple rank-1 lattices.

Suitable function spaces $\mathcal{F}$ and $Y$ allow for detailed estimates
of the corresponding sampling numbers. In this paper we mainly deal with commonly used
Hilbert spaces, cf.\ e.g.\ \cite{ByDuSiUl14, ByKaUlVo16}.
Using the notation from~\cite{KaPoVo13},
we consider the periodic Sobolev spaces of generalized mixed smoothness %
\begin{equation*} %
 \mathcal{H}^{\alpha,\beta}(\T^d) := \left\{ f\in L_1(\T^d)\colon \Vert f\vert \mathcal{H}^{\alpha,\beta}(\T^d) \Vert := \sqrt{ \sum_{\boldk\in\Z^d} \omega^{\alpha,\beta}(\boldk)^2 \vert\hat{f}_\boldk\vert^2 } < \infty \right\}
\end{equation*}
$ \subset L_2(\T^d)$ with dominating mixed smoothness $\beta\geq 0$, isotropic smoothness $\alpha \geq -\beta$,
where the weights $\omega^{\alpha,\beta}\colon\Z^d\to(0,\infty)$ are defined by
\begin{equation} \label{equ:weight:omega}
 \omega^{\alpha,\beta}(\boldk) := \max(1,\Vert\boldk\Vert_1)^\alpha \; \prod_{s=1}^{d} \max(1,\vert k_s\vert)^\beta, \quad
 \boldk := \left(\substack{k_1\\\vdots\\k_d}\right),
\end{equation}
and the Fourier coefficients~$\hat{f}_\boldk$ of~$f$
are formally given by~\eqref{equ:fc:integral}.

We remark that these function spaces $\mathcal{H}^{\alpha,\beta}(\T^d)$ are Hilbert spaces
and that $\mathcal{H}^{0,0}(\T^d)$ coincides with the Lebesgue space $L_2(\T^d)$ of all square integrable functions.
Moreover, we denote the Hilbert spaces of Sobolev type containing functions of specific isotropic smoothness by $\mathcal{H}^{r}(\T^d):=\mathcal{H}^{r,0}(\T^d)$
as well as the Hilbert spaces of Sobolev type containing functions of specific dominating mixed smoothness by $\mathcal{H}^{t}_{\mathrm{mix}}(\T^d):=\mathcal{H}^{0,t}(\T^d)$.
Note that for positive integers $r$ the spaces $\mathcal{H}^{r}(\T^d):=\mathcal{H}^{r,0}(\T^d)$ consist
of all functions whose weak derivatives up to total degree $r$ are contained in $L_2(\T^d)$.
Furthermore, functions with weak derivatives in $L_2(\T^d)$ up to degree $t\in\N$ in each variable constitute the spaces
$\mathcal{H}^{t}_{\mathrm{mix}}(\T^d):=\mathcal{H}^{0,t}(\T^d)$.

The second columns of Tables~\ref{tab:mainres:overview_general_r1l} and~\ref{tab:mainres:overview_general_r1l_2}
present the main contributions of this paper,
i.e., upper bounds on sampling numbers $g_M^{\mathrm{mr1l}}$ for multiple \hbox{rank-1} lattices for important parameter combinations.
Additionally, these results  are compared with known bounds on the sampling numbers for
single rank-1 lattices $g_M^{\mathrm{latt}_1}$ and on the sampling numbers for linear operators $g_M^{\mathrm{lin}}$, cf.\ the third and fourth columns of Tables~\ref{tab:mainres:overview_general_r1l} and~\ref{tab:mainres:overview_general_r1l_2}.
The important aspect for the comparisons of the bounds is the exponent of~$M$, called main rate.
For Hilbert spaces $Y\in\{L_2(\T^d), \mathcal{H}^{r}(\T^d), \mathcal{H}^{t}_{\mathrm{mix}}(\T^d)\}$ as target space,
the main rates of the newly introduced sampling numbers $g_M^{\mathrm{mr1l}}$ for multiple rank-1 lattices
correspond to the best known bounds on the sampling numbers for linear operators $g_M^{\mathrm{lin}}$
up to an offset slightly larger than one half.
Compared to the sampling numbers $g_M^{\mathrm{latt}_1}$ for single rank-1 lattices,
which only yield half of the main rates in lower and upper bounds, cf.~\cite{ByKaUlVo16},
this means a distinct improvement for lattice based sampling.

Considering the target space $Y=L_\infty(\T^d)$, the main rate of
the sampling numbers $g_M^{\mathrm{mr1l}}$ for multiple rank-1 lattices
is optimal up to an arbitrarily small offset $\varepsilon>0$,
whereas the main rates of the sampling numbers $g_M^{\mathrm{latt}_1}$ for single rank-1 lattices
are only half of the optimal ones for linear operators, cf.\ the third rows in
Tables~\ref{tab:mainres:overview_general_r1l} and~\ref{tab:mainres:overview_general_r1l_2}.

\begin{table}[ht]
	\begin{center}
	\scalebox{0.95}{%
		\begin{tabular}{l|lll}
			\toprule
			$Y$ & $g_M^{\mathrm{mr1l}}(\mathcal{H}^{\beta}_{\text{mix}}(\T^d),Y)$ & $g_M^{\mathrm{latt}_1}(\mathcal{H}^{\beta}_{\text{mix}}(\T^d),Y)$ &
			$g_M^{\mathrm{lin}}(\mathcal{H}^{\beta}_{\text{mix}}(\T^d),Y)$ \\
			& \scriptsize $\left(\begin{array}{l}\beta-\frac{1}{2}-\varepsilon>\max\{r,t\} \ge 0,\\ \textnormal{and } \varepsilon>0\end{array}\right)$  &&\\
			\cmidrule{1-4}
			$L_2(\T^d)$  &$\lesssim 
			\left(\frac{\log^d M}{M}\right)^{\beta-\frac{1}{2}-\varepsilon}
			$
			& $\lesssim 
			\left(\frac{\log^{d-2} M}{M}\right)^{\frac{\beta}{2}}
			$ & $\lesssim 
			\left(\frac{\log^{d-1} M}{M}\right)^\beta 
			$ \\
			\rule[0em]{0em}{1.2em} & $\quad\;\cdot \log M$ & $\quad\;\cdot (\log M)^{\frac{d-1}{2}}$ & $\quad\;\cdot (\log M)^{\frac{d-1}{2}}$ \\
			&
			{\footnotesize (Corollary~\ref{cor:approx_error:Hrtg:Habg})}
			&{\footnotesize \cite[Theorem 2]{ByKaUlVo16}}&{
				{\footnotesize\cite[Theorem 6.10]{ByDuSiUl14}}}\\
			\cmidrule{1-4}
			$L_{\infty}(\T^d)$  &$\lesssim 
			\left(\frac{\log^d M}{M}\right)^{\beta-\frac{1}{2}-\varepsilon}
			$
			& $\lesssim
			\left(\frac{\log^{d-2} M}{M}\right)^{\frac{\beta}{2}-\frac{1}{4}}
			$
			& $\asymp
			\left(\frac{\log^{d-1} M}{M}\right)^\beta M^{\frac{1}{2}}
			$ \\
			\rule[0em]{0em}{1.2em} & $\quad\;\cdot \log M$ & $\quad\;\cdot (\log M)^{\frac{d-1}{2}}$ & \\
			&
			{\footnotesize (Corollary~\ref{cor:sampl_error_mr1l:Linf:Aabg} \& Lemma~\ref{lemma:estimate:Aabg:Habg})}
			&{\footnotesize \cite[Proposition 4]{ByKaUlVo16}}&{
				{\footnotesize\cite[Theorem 6.10]{ByDuSiUl14}}}\\
			\cmidrule{1-4}
			$\mathcal{H}^{r}(\T^d)$&  $\lesssim
			\left(\frac{\log^d M}{M}\right)^{\beta-r-\frac{1}{2}-\varepsilon}
			$			
			&$\lesssim 
			\left(\frac{\log^{d-2} M}{M}\right)^{\frac{\beta-r}{2}}
			$ & $\asymp
			M^{-(\beta-r)}$\\
			\rule[0em]{0em}{1.2em} & $\quad\;\cdot \log M$ & $\quad\;\cdot (\log M)^{\frac{d-1}{2}}$ & \\
			&
			{\footnotesize (Corollary~\ref{cor:approx_error:Hrtg:Habg})}
			&{\footnotesize \cite[Proposition 2]{ByKaUlVo16}}&{
				{\footnotesize\cite[Theorem 6.7]{ByDuSiUl14}}}\\
			\cmidrule{1-4}
			$\mathcal{H}^{t}_{\mathrm{mix}}(\T^d)$& $\lesssim
						\left(\frac{\log^d M}{M}\right)^{\beta-t-\frac{1}{2}-\varepsilon}
						$			
			& $\lesssim
			\left(\frac{\log^{d-2} M}{M}\right)^{\frac{\beta-t}{2}}
			$ & $\asymp
			\left(\frac{\log^{d-1} M}{M}\right)^{\beta-t}$
			\\
			\rule[0em]{0em}{1.2em} & $\quad\;\cdot \log M$ & $\quad\;\cdot (\log M)^{\frac{d-1}{2}}$ & \\
			&
			{\footnotesize (Corollary~\ref{cor:approx_error:Hrtg:Habg})}
			&{\footnotesize \cite[Theorem 2]{ByKaUlVo16}}&{
				{\footnotesize\cite[Theorem 6.10]{ByDuSiUl14}}}\\
			\bottomrule
		\end{tabular}
}
	\end{center}
	\caption{Upper bounds of sampling numbers in the setting
		$\mathcal{H}^{\beta}_{\text{mix}}(\T^d)\rightarrow Y$ for different sampling
		methods.
		Smoothness parameters are chosen from
		$\beta>\max\{r,t,\frac{1}{2}\}$ and $r,t>0$.
		}\label{tab:mainres:overview_general_r1l}
\end{table}
\begin{table}[ht]
	\centering
	\scalebox{0.93}{%
	\begin{tabular}{l|lll}
		\toprule
		Y & $g_M^{\mathrm{mr1l}}(\mathcal{H}^{\alpha,\beta}(\T^d),Y)$ & $g_M^{\mathrm{latt}_1}(\mathcal{H}^{\alpha,\beta}(\T^d),Y)$ &
		$g_M^{\mathrm{lin}}(\mathcal{H}^{\alpha,\beta}(\T^d),Y)$ \\
		& \scriptsize $\left(\begin{array}{l}\alpha+\beta-\frac{1}{2}-\varepsilon>\max\{r,t\} \ge 0,\\ \textnormal{and } \varepsilon>0\end{array}\right)$  &&\\
		\cmidrule{1-4}
		$ L_2(\T^d)$  &
		$\lesssim 
		\left(\frac{\log M}{M}\right)^{\alpha+\beta-\frac{1}{2}-\varepsilon} \log M
		$		
		& $\lesssim M^{-\frac{\alpha+\beta}{2}}$ & $\asymp M^{-(\alpha+\beta)}$ \\
		&
		{\footnotesize (Corollary~\ref{cor:approx_error:Hrtg:Habg})}
		&{\footnotesize \cite[Theorem 4.7]{KaPoVo14}}&{ {\footnotesize\cite[Theorem 6.10]{ByDuSiUl14}}}\\
		\cmidrule{1-4}
		$L_{\infty}(\T^d)$  &
		$\lesssim 
		\left(\frac{\log M}{M}\right)^{\alpha+\beta-\frac{1}{2}-\varepsilon} \log M
		$				
		& $\lesssim M^{-\frac{\alpha+\beta-\frac{1}{2}}{2}}$ & $\lesssim M^{-(\alpha+\beta)+\frac{1}{2}}$ \\
		&
		{\footnotesize (Corollary~\ref{cor:sampl_error_mr1l:Linf:Aabg} and Lemma~\ref{lemma:estimate:Aabg:Habg})}
		&{\footnotesize \cite[Proposition 4]{ByKaUlVo16}}&
		{ {\footnotesize\cite[*]{ByDuSiUl14}}}\\
		\cmidrule{1-4}
		$\mathcal{H}^{r}(\T^d)$& 
		$\lesssim 
		\left(\frac{\log M}{M}\right)^{\alpha-r+\beta-\frac{1}{2}-\varepsilon} \log M
		$				
		&  $\lesssim 
		\left(\frac{\log^{d-2}M}{M}\right)^{\frac{\alpha-r+\beta}{2}}
		$ 
		& $\asymp
		M^{-(\alpha-r+\beta)}$ \\
		&
		{\footnotesize (Corollary~\ref{cor:approx_error:Hrtg:Habg})}
		&{\footnotesize \cite[Proposition 2]{ByKaUlVo16}}&{ {\footnotesize\cite[Theorem 6.7]{ByDuSiUl14}}}\\
		\cmidrule{1-4}
		$\mathcal{H}^{t}_{\mathrm{mix}}(\T^d)$& 
		$\lesssim 
		\left(\frac{\log M}{M}\right)^{\alpha+\beta-t-\frac{1}{2}-\varepsilon} \log M
		$		
		&  $\lesssim
		\left(\frac{\log^{d-2}M}{M}\right)^{\frac{\alpha+\beta-t}{2}}
		$ & $\asymp
		M^{-(\alpha+\beta-t)}$ \\
		&
		{\footnotesize (Corollary~\ref{cor:approx_error:Hrtg:Habg})}
		&{\footnotesize \cite[Theorem 2]{ByKaUlVo16}}&{ {\footnotesize\cite[*]{ByDuSiUl14}}}\\
		\bottomrule
	\end{tabular}
}
	\caption{Upper bounds for sampling numbers for different sampling methods.
	Smoothness parameters are chosen from
	$\alpha<0$, $\alpha+\beta>\max\{r,t,\frac{1}{2}\}$ and $r,t>0$. Best known bounds based on energy sparse grid sampling.
	References marked with $^*$ mean that the result is not stated there explicitly but follows with the same method
	therein.}
	\label{tab:mainres:overview_general_r1l_2}
	\end{table}

The remaining parts of this paper are organized as follows.
In Section~\ref{sec:pre}, we collect important facts on used function spaces and frequency index sets.
Moreover, we recapitulate important, already known facts on the reconstruction of arbitrary multivariate trigonometric polynomials based on samples along single and multiple rank-1 lattices.
Subsequently, we prove upper bounds on worst case errors for the presented sampling strategy based on reconstructing multiple rank-1 lattices in Section~\ref{sec:mult_r1l_sampl}
and these yield the upper bounds on sampling numbers~$g_M^{\ml}$ for sampling sets consisting of multiple rank-1 lattices as presented in Tables~\ref{tab:mainres:overview_general_r1l} and~\ref{tab:mainres:overview_general_r1l_2}. In Section~\ref{sec:numerics}, we verify our theoretical results in numerical tests and compare the numerical results with those for single rank-1 lattice sampling and sparse grid sampling.

\section{Prerequisites}\label{sec:pre}

\subsection{Function spaces and frequency index sets}

For our theoretical considerations, we introduce the subspaces 
\begin{equation*} %
 \mathcal{A}^{\alpha,\beta}(\T^d) := \left\{ f\in L_1(\T^d)\colon \Vert f\vert \mathcal{A}^{\alpha,\beta}(\T^d) \Vert := \sum_{\boldk\in\Z^d} \omega^{\alpha,\beta}(\boldk) \vert\hat{f}_\boldk\vert < \infty \right\}
\end{equation*}
of the Wiener algebra $\mathcal{A}(\T^d)=\mathcal{A}^{0,0}(\T^d)$
with dominating mixed smoothness $\beta\geq 0$, isotropic smoothness $\alpha \geq -\beta$,
where the weights $\omega^{\alpha,\beta}$ are defined in~\eqref{equ:weight:omega}.
We remark that the embedding $\mathcal{A}(\T^d)\hookrightarrow\mathcal{C}(\T^d)$ holds, cf.\ e.g.\ \cite[Remark~3.1]{volkmerdiss}.

In Section~\ref{sec:mult_r1l_sampl}, we require the following embeddings between subspaces $\mathcal{A}^{\alpha,\beta}(\T^d)$ of the Wiener algebra $\mathcal{A}(\T^d)$ and
periodic Sobolev spaces of generalized mixed smoothness~$\mathcal{H}^{\alpha,\beta}(\T^d)$.
\begin{Lemma} \label{lemma:estimate:Aabg:Habg}
 (\cite[Lemma 2.2]{KaPoVo13}).
 Let a function $f\in\mathcal{A}^{\alpha,\beta}(\T^d)$ be given,
 where
 the dominating mixed smoothness $\beta\geq 0$
 and the isotropic smoothness
 $\alpha > -\beta$.
 Then, we have
 $ \Vert f\vert \mathcal{H}^{\alpha,\beta}(\T^d) \Vert
   \leq
   \Vert f\vert \mathcal{A}^{\alpha,\beta}(\T^d) \Vert.
 $
 For $f\in\mathcal{H}^{\alpha,\beta+\lambda}(\T^d)$ with $\lambda>1/2$,
 we have
 \begin{equation*} %
  \Vert f\vert \mathcal{A}^{\alpha,\beta}(\T^d) \Vert
  \leq
	\prod_{s=1}^d( 1+2\zeta(2\lambda))^{1/2}
  \Vert f\vert \mathcal{H}^{\alpha,\beta+\lambda}(\T^d) \Vert,
 \end{equation*}
where $\zeta$ denotes the Riemann zeta function.
 Therefore, we have the continuous embeddings
 $$
 \mathcal{H}^{\alpha,\beta+\lambda}(\T^d)
 \hookrightarrow
 \mathcal{A}^{\alpha,\beta}(\T^d)
 \hookrightarrow
 \mathcal{A}(\T^d)
 \hookrightarrow
 L_2(\T^d)
 \hookrightarrow
 L_1(\T^d).
 $$
\end{Lemma}

If the Fourier coefficients~$\hat{f}_\boldk$ decay in a certain way, one may approximate~$f$ by a Fourier partial sum $S_I f$
with respect to a frequency index set $I\subset\Z^d$, $|I|<\infty$, defined by
$$S_If:=\sum_{\zb k\in I}\hat{f}_{\zb k}\e^{2\pi\ii\zb k\cdot\circ}.$$
As frequency index sets $I$,
we use %
 \begin{equation} \label{equ:I_N_d_T}
  I_N^{d,T} := \left\{\boldk\in\Z^d\colon \omega^{-T,1}(\boldk) = \max(1,\Vert\boldk\Vert_1)^{-T} \; \prod_{s=1}^{d} \max(1,\vert k_s\vert) \leq N^{1-T} \right\},
 \end{equation}
where
$N\geq 1$ is the refinement, $T\in (-\infty,1)$ is the shape parameter,
and the weights $\omega^{\alpha,\beta}(\boldk)$ are defined as in~\eqref{equ:weight:omega}.
As a natural extension for $T=-\infty$, we define the frequency index set $I_N^{d,-\infty}$
as the $d$-dimensional $\ell_1$-ball of size~$N$,
 \begin{equation*} %
  I_N^{d,-\infty} := \left\{\boldk\in\Z^d\colon \max(1,\Vert\boldk\Vert_1)\leq N \right\}.
 \end{equation*}
The number of degrees of freedom when reconstructing $S_If$
is $|I|$. Correspondingly, we bound
the cardinalities of the frequency index sets $I_N^{d,T}$.

\begin{Lemma} \label{lemma:index_set:omega}
 (\cite[Lemma~4.1]{KaPoVo14}).
 Let the dimension $d\in\N$, and a shape parameter $T\in [-\infty,1)$ be given.
 Then, the cardinalities of the frequency index sets $I_N^{d,T}$ are
\begin{equation*} %
 \vert I_N^{d,T} \vert
 \asymp
 \begin{cases}
  N^d & \textrm{for } T = -\infty,\\
  N^{\frac{T-1}{T/d-1}} & \textrm{for } -\infty < T < 0, \\
  N \log^{d-1} N & \textrm{for } T = 0, \\
  N & \textrm{for } 0 < T < 1, \\
 \end{cases}
\end{equation*}
 for fixed parameters $d$ and $T$, where the constants only depend on $d$ and $T$.
\end{Lemma}

However, in practice, the Fourier coefficients~$\hat{f}_\boldk$ of a function~$f$ are often hard to compute or even unknown.
Then, one may approximate the Fourier coefficients~$\hat{f}_\boldk$ based on sampling values of~$f$.
Two possible sampling strategies are considered in the next sections.

\subsection{Reconstructing (single) rank-1 lattices}
\label{sec:pre:single_r1l}

First, we start with definitions from \cite{Kae17} using slightly adapted symbols in this work.

The sampling sets $\mathcal{G}$ which are used in this work are based on so-called
rank-1 lattices $\Lambda(\boldz,M)$ as defined in~\eqref{eq:def:r1l},
where
$\boldz\in\Z^d$ and $M\in\N$ are called generating vector and lattice size of $\Lambda(\boldz,M)$, respectively.
For an arbitrary multivariate trigonometric polynomial
\[
p\in\Pi_I:=\sspan\{\e^{2\pi\ii\zb k\cdot\circ}\colon\zb k\in I\}, \quad
p(\boldx):= \sum_{\boldk\in I} \hat{p}_\boldk \,\mathrm{e}^{2\pi\mathrm{i}\boldk\cdot\boldx}, \; \hat{p}_\boldk\in\C,
\]
with frequencies~$\boldk$ supported on an index set $I\subset\Z^d$, $|I|<\infty$, we can reconstruct all the Fourier coefficients $\hat{p}_\boldk$, $\boldk\in I$, %
from samples along a \hbox{rank-1} lattice $\mathcal{G}:=\Lambda(\boldz,M)$ if the Fourier matrix $
\boldsymbol{A}(\mathcal{G},I):=\left( 
\mathrm{e}^{2\pi\mathrm{i}\boldk\cdot\boldx} \right)_{\boldx\in \mathcal{G}, \boldk \in I}
$
has full column rank, see e.g.\ \cite{Kae2013,kaemmererdiss}.
This is the case if and only if $\Lambda(\boldz,M)$ is a \textit{reconstructing rank-1 lattice} for $I$, i.e., the \textit{reconstruction property}
\begin{equation}\label{eqn:reco_prop1}
 \boldk\cdot\boldz\not\equiv \boldsymbol{k'}\cdot\boldz\imod{M} \text{ for all } \boldk,\boldsymbol{k'}\in I, \boldk\neq\boldsymbol{k'},
\end{equation}
is fulfilled.

Using sampling values of $p$ along a reconstructing rank-1 lattice~$\Lambda(\boldz,M)$ for~$I$, the reconstruction can be performed in a fast way by applying a one-dimensional fast Fourier transform,
\begin{align}
\hat{h}_l:=&\sum_{j=0}^{M-1} p(\frac{j}{M}\boldz\bmod{\boldsymbol{1}}) \; \mathrm{e}^{-2\pi\mathrm{i}l j/M},\qquad l=0,\ldots,M-1,\label{eq:r1lfft_1}\\
\intertext{followed by the simple index transform}
\hat{p}^{\,\Lambda(\boldz,M)}_\boldk:=&\frac{1}{M} \, \hat{h}_{\boldk\cdot\boldz\bmod{M}}, \qquad \boldk\in I,\label{eq:r1lfft_2}
\end{align}
cf.\ \cite[Algorithm~3.2]{kaemmererdiss}.
This computation requires $\mathcal{O}(M \log M + d\,\vert I\vert)$
arithmetic operations.

Moreover, a reconstructing rank-1 lattice~$\Lambda(\boldz,M)$ can be easily constructed using a simple component-by-component construction method, cf.~\cite{Kae2013}.
However, the construction method has rather high computational costs and may require $\mathcal{O}(d|I|^3)$ arithmetic operations.
Furthermore, under mild assumptions, the number of samples $M$ is bounded by $|I|\leq M \leq |I|^2$, where this number tends more to the upper bound for many interesting structures of frequency index sets~$I$.
For instance, for the axis cross $I:=\{\boldk\in\Z^d\colon \|\boldk\|_\infty=\|\boldk\|_1 \leq N \}$, $N\in\N$, it can be shown that $M\gtrsim |I|^2$ is a necessary condition for a reconstructing rank-1 lattice, cf.~\cite{kaemmererdiss,ByKaUlVo16}.

Applying \eqref{eq:r1lfft_1} and \eqref{eq:r1lfft_2} to sampling values of continuous periodic functions~$f\colon\T^d\rightarrow\C$ for a given frequency index set $I\subset\Z^d$ and rank-1 lattice $\Lambda(\boldz,M)$,
we obtain all the approximated Fourier coefficients
\begin{align} \label{equ:fhattilde}
 \hat{f}^{\,\Lambda(\boldz,M)}_\boldk :=& \frac{1}{M} \sum_{j=0}^{M-1} f\left(\frac{j}{M}\boldz\bmod{\boldsymbol{1}}\right) \, \mathrm{e}^{-2\pi\mathrm{i}j\boldk\cdot\boldz/M} \\
 =& \sum_{\substack{\boldsymbol{k'}\in\Z^d \\ \boldk\cdot\boldz \equiv \boldsymbol{k'}\cdot\boldz \bmod{M}}} \hat{f}_{\boldsymbol{k'}} = \sum_{\boldh\in\Lambda(\zb z,M)^\perp} \hat{f}_{\boldk+\boldh}
 \qquad \boldk\in I, \nonumber
\end{align}
where
\begin{equation}
\Lambda(\zb z,M)^\perp:=\{\zb h\in\Z^d\colon\zb h\cdot\zb z\equiv 0\imod{M}\}
\label{equ:def_dual_lattice}
\end{equation}
is the (integer) dual lattice 
of the rank-1 lattice $\Lambda(\zb z,M)$.
Based on the approximated Fourier coefficients~$\hat{f}^{\,\Lambda(\boldz,M)}_\boldk$,
we define the 
rank-1 lattice sampling operator $S_I^{\Lambda(\zb z,M)}$ by
\begin{align}
S_I^{\Lambda(\zb z,M)}f&:=\sum_{\zb k\in I}\hat{f}_{\zb k}^{\Lambda(\zb z,M)}\e^{2\pi\ii\zb k\cdot\circ}\label{eqn:sampling_operator}.
\end{align}
We stress on the fact that the reconstruction property~\eqref{eqn:reco_prop1} is equivalent to the requirement that the sampling operator $S_I^{\Lambda(\zb z,M)}$ reproduces all trigonometric polynomials~$p$ with frequencies supported on I, i.e.,
$S_I^{\Lambda(\zb z,M)}p=p$ holds for
$p\in\Pi_I$. %

For various approximation settings, the errors for single rank-1 lattice sampling have already been investigated, cf.\ e.g.\ \cite{Tem86,KuSlWo06,KuSlWo08,KaPoVo13,KaPoVo14,volkmerdiss,ByKaUlVo16}.
When estimating the sampling error $f - S^{\Lambda(\boldz,M)}_I f$,
one usually splits this error into the truncation error and aliasing error,
$
f - S^{\Lambda(\boldz,M)}_I f = (f - S_I f) + (S_I f - S^{\Lambda(\boldz,M)}_I f).
$
Applying the triangle inequality on the norm $\Vert \circ | Y \Vert$ of the target space $Y$ yields
\begin{equation}
\label{equ:periodic:approx:triangle}
\Vert f - S^{\Lambda(\boldz,M)}_I f |Y\Vert \leq \Vert f - S_I f |Y\Vert + \Vert S_I f - S^{\Lambda(\boldz,M)}_I f |Y\Vert.
\end{equation}
If the employed single rank-1 lattice $\Lambda(\zb z,M)$ is a reconstructing one for the frequency index set~$I$,
theoretical upper bounds on the truncation error $\Vert f - S_I f |Y\Vert$ and aliasing error $\Vert S_I f - S^{\Lambda(\boldz,M)}_I f |Y\Vert$ are of comparable order of magnitude in many cases.
However, the main issue when using reconstructing single rank-1 lattices $\Lambda(\boldz,M)$ as sampling schemes is the (asymptotically) large oversampling $M \gg |I|$
for the arising structures of frequency index sets~$I$ in highly interesting approximation settings.
This large oversampling leads to sampling errors of lower order with respect to $M$, compare e.g.\ the sampling numbers~$g_M^{\mathrm{latt}_1}$ for single rank-1 lattice sampling with the sampling numbers~$g_M^{\mathrm{lin}}$ for general linear sampling operators in Tables~\ref{tab:mainres:overview_general_r1l} and~\ref{tab:mainres:overview_general_r1l_2}.
We remark that the lower bounds on the sampling numbers~$g_M^{\mathrm{latt}_1}$
correspond to the upper bounds in the main rate, \cite[Section~3]{ByKaUlVo16}.

\subsection{Reconstructing multiple rank-1 lattices}
\label{sec:pre:multiple_r1l}

Recently, in \cite{Kae16,Kae17}, a modified approach was presented, which allows one to drastically reduce the number of samples when reconstructing arbitrary multivariate trigonometric polynomials~$p$. This approach uses rank-1 lattices $\Lambda(\boldz,M)$ as building blocks and builds sampling sets based on multiple instances. The corresponding sampling sets $\mathcal{G}$ are called multiple rank-1 lattices and they can be constructed by simple and efficient randomized construction algorithms, cf.\ \cite[Algorithms~1 to~6]{Kae17}.
A multiple rank-1 lattice $\Lambda=\Lambda(\boldz_1,M_1,\ldots,\boldz_L,M_L)$ is the union of $L\in\N$ single rank-1 lattices, cf.~\eqref{eq:def:mr1l},
and consists of
$
|\Lambda(\boldz_1,M_1,\ldots,\boldz_L,M_L)|
\leq
1-L+\sum_{\ell=1}^{L}M_\ell
$
distinct nodes.
If~$\Lambda$ allows for the reconstruction of all multivariate trigonometric polynomials $p$ with frequencies supported on a frequency index set $I$,
it will be called \textit{reconstructing multiple rank-1 lattice} for $I$.

In simplified terms, the basic idea is that each of the rank-1 lattices $\Lambda(\boldz_\ell,M_\ell)$, $\ell=1,\ldots,L$, should be a reconstructing one for some index set $I_\ell\subset I$
and that $\bigcup_{\ell=1}^L I_\ell = I$. We remark that this condition is not sufficient in general and we require additional properties.

\begin{algorithm}[tb]
\caption{(\cite[Algorithm~4]{Kae17}). Determining reconstructing multiple rank-1 lattices  with pairwise distinct lattice sizes $M_\ell$
}\label{alg:construct_mr1l_I_distinct_primes}
  \begin{tabular}{p{2.2cm}p{4.0cm}p{7.cm}}
    Input: 	& $I\subset\Z^d$ 	& frequency index set\\
    		& $c\in(1,\infty)\subset\R$ 	& oversampling factor\\
    		& $\delta\in(0,1)\subset\R$			& upper bound on failure probability
  \end{tabular}
		
  \begin{algorithmic}
  	\STATE $L_{\mathrm{max}}:=\ceil{\left(\frac{c}{c-1}\right)^2\frac{\ln T-\ln\delta}{2}}$
  	\STATE $\lambda:=c(T-1)$
	\STATE determine $P^{I}_{\lambda,L_{\mathrm{max}}}$, cf.~\eqref{eq:def_PIlambdan}, and arrange $p_1<\ldots<p_{L_{\mathrm{max}}}$
	\STATE $\tilde{I}:=\emptyset$
	\STATE $L:=0$
	\WHILE {$|\tilde{I}|<|I|$ and $L<L_{\mathrm{max}}$}
		\STATE $L:=L+1$
		\STATE choose $M_{L}:=p_{L}\in P^{I}_{\lambda,L_{\mathrm{max}}}$\label{alg:construct_mr1l_I_distinct_primes:determinstic_not_random_lattice_sizes}
		\STATE choose $\zb z_{L}$ from $[0,M_{L}-1]^d\cap\Z^d$ uniformly at random
		\IF {$\{\zb k\in I\colon \not\exists \zb h\in I\setminus\{\zb k\}\text{ with }\zb k\cdot \zb z_{L}\equiv\zb h\cdot \zb z_{L}\imod{M_{L}}\}\not\subset \tilde{I}$}
	    \STATE compute \\
	    {\hfill\mbox{$\tilde{I}:=\tilde{I}\cup \{\zb k\in I\colon \not\exists \zb h\in I\setminus\{\zb k\}\text{ with }\zb k\cdot \zb z_{L}\equiv\zb h\cdot \zb z_{L}\imod{M_{L}}\}$}}
	    \ELSE\label{alg:construct_mr1l_I_distinct_primes_line:avoid_useless_r1l_1}
			\STATE	$L:=L-1$
	    \ENDIF\label{alg:construct_mr1l_I_distinct_primes_line:avoid_useless_r1l_2}
	\ENDWHILE
  \end{algorithmic}
  \begin{tabular}{p{2.2cm}p{4.2cm}p{8.525cm}}
    Output: & $M_1,\ldots,M_{L}$ & lattice sizes of rank\mbox{-}1 lattices and\\
	    & $\zb z_1,\ldots,\zb z_{L}$ & generating vectors of rank\mbox{-}1 lattices such that\\
	    & $\Lambda(\zb z_1,M_1,\ldots,\zb z_{L},M_{L})$ &  is a reconstructing multiple rank\mbox{-}1 lattice
														   for $I$ with probability at least $1-\delta$  
\\    \cmidrule{1-3}
  \end{tabular}
  \begin{tabular}{p{2.2cm}p{11.5cm}}
    Complexity: & $\OO{|I|(d+\log|I|)\log|I|}$\; w.h.p.\\& {for $c(|I|-1)\geq N_I$ and fixed $c$ and $\delta$, where\newline $N_I:=\max_{j=1,\ldots, d}\{\max_{\zb k\in I}k_j-\min_{\zb l\in I}l_j\}$ is the expansion of~$I$}
  \end{tabular}
\end{algorithm}

The construction approach utilized in this work is Algorithm~\ref{alg:construct_mr1l_I_distinct_primes}, originally stated as \cite[Algorithm~4]{Kae17}, which determines a reconstructing multiple rank-1 lattice $\Lambda$ for a given index set $I$ with high probability,
such that the reconstruction property
\begin{equation}\label{equ:cond:reco_ml_mlfft2}
\bigcup_{\ell=1}^L I_\ell = I
\end{equation}
is fulfilled with index sets %
$$I_\ell:=\left\{\boldk \in I\colon \boldk\cdot\boldz_\ell\not\equiv\boldk'\cdot\boldz_\ell\imod{M_\ell}\textnormal{ for all }\boldk'\in I\setminus\{\boldk\}\right\}.$$
The rank-1 lattice sizes $M_\ell$ are chosen distinctly from the set
\begin{align}
P^{I}_{\lambda,L_\mathrm{max}}&:=\Bigg\{p_j\in P^I\colon p_j:=
\text{\scalebox{0.85}{
$\displaystyle\begin{cases}
\min\{p \in P^I\colon p>\lambda\}&\colon j=1\\
\min\{p \in P^I\colon p>p_{j-1}\}&\colon j=2,\ldots,L_\mathrm{max}.
\end{cases}$}}
 \Bigg\}\label{eq:def_PIlambdan}
\end{align}
of the $L_\mathrm{max}\in\N$ smallest prime numbers in
$$
P^I:=\{M'\in\N\colon M' \text{ prime with }|\{\zb k\bmod M'\colon \zb k\in I\}|=|I|\}
$$
larger than a certain $\lambda\in\N$.
Under mild technical assumptions, cf.\ \cite[Corollary~3.7]{Kae17}, the mentioned algorithm returns a multiple rank-1 lattice $\Lambda=\Lambda(\zb z_1,M_1,\ldots,\zb z_L,M_L)$ of cardinality $|\Lambda|\le M \in\mathcal{O}(|I|\log|I|)$, where $M:=\sum_{\ell=1}^{L}M_\ell$ is an upper bound on the number of sampling nodes within $\Lambda$, the lattice sizes $M_\ell\approx c\,|I|$, $c>1$, and $L\lesssim\log|I|$.
With high probability, $\Lambda$ fulfills reconstruction property~\eqref{equ:cond:reco_ml_mlfft2} and the construction requires $\mathcal{O}(|I| (d+\log|I|) \log|I|)$ arithmetic operations, see also \cite[Corollary~3.7]{Kae17}.
At this point, we stress the facts that the oversampling factor $\frac{1+\sum_{\ell=1}^L (M_{\ell}-1)}{|I|}\le\frac{M}{|I|}\lesssim \log|I|$ does not depend on the dimension~$d$ and that
checking the reconstruction property~\eqref{equ:cond:reco_ml_mlfft2} can be efficiently
performed during construction in Algorithm~\ref{alg:construct_mr1l_I_distinct_primes}.

Besides the small cardinalities and fast construction algorithms, a further main advantage of reconstructing multiple rank-1 lattices is the existence of a direct and fast inversion method for computing Fourier coefficients $\hat{p}_\boldk$ from samples, cf.\ Algorithm~\ref{algo:mlfft2:ifft_direct}.
Its arithmetic complexity is $\mathcal{O}\big(M\log M + L\,|I|\,(d + \log |I|)\big)$, 
where the second summand comes from the computation of the index sets $I_\ell$ in line~\ref{alg:mlfftinv_Il_line}
and the complexity of the adjoint rank-1 lattice  FFTs in line~\ref{alg:mlfftinv_lfft_line}
is $\mathcal{O}\big(M\log M + L\,|I|\,d\big)$ in total.
For more details on the adjoint rank-1 lattice FFT, we refer to \cite[Algorithm~3.2]{kaemmererdiss}.
The computation of the index sets $I_\ell$ can be performed similarly as in the construction of
reconstructing multiple rank-1 lattices, cf.\ \cite[Section~3]{Kae17} in the context of \cite[Equation~3.8]{Kae17} for full details.

\begin{algorithm}[t]
\caption{(\cite[Algorithm~2]{KaPoVo17}). Reconstruction of a multivariate trigonometric polynomial $p$ from sampling values along reconstructing multiple rank-1 lattice fulfilling reconstruction property~\eqref{equ:cond:reco_ml_mlfft2}.}\label{algo:mlfft2:ifft_direct}
  \scalebox{0.98}{%
  \begin{tabular}{p{2.2cm}p{5.4cm}p{7.6cm}}
    Input:      & $I\subset\Z^d$ & frequency index set, $|I|<\infty$, \\
                & $\Lambda:=\Lambda(\boldz_1,M_1,\ldots,\boldz_L,M_L)$ & reconstructing multiple rank-1 lattice for $I$ obtained from one of \cite[Algorithms~1--4]{Kae17}, \\
                & $\left(p(\boldsymbol{\tilde{x}}_j)\right)_{\boldsymbol{\tilde{x}}_j\in\Lambda}$  & samples of trigonometric polynomial $p$ \\
  \end{tabular}}
  \begin{algorithmic}[1]
     \STATE Initialize $\mathtt{counter}[\boldk]:=0$ and\newline set output Fourier coefficients $\hat{p}^{\,\Lambda}_\boldk:=0$ for $\boldk\in I$.
	\FORALL{$\ell\in\{1,\ldots,L\}$}
     \STATE Determine frequency index set\newline
     {\centering
     $I_\ell:=\{\zb k\in I\colon \nexists \zb h\in I\setminus\{\zb k\}\text{ with }\zb k\cdot \zb z_\ell\equiv\zb h\cdot \zb z_\ell\imod{M_\ell}\}, \textnormal{ cf.~\eqref{equ:cond:reco_ml_mlfft2}},$}\newline
     containing those frequencies~$\boldk\in I$ where the corresponding Fourier coefficients~$\hat{p}_\boldk$ of $p\in\Pi_I$ can be exactly reconstructed using the samples $\left(p(\boldsymbol{\tilde{x}}_j)\right)_{\boldsymbol{\tilde{x}}_j\in\Lambda(\boldz_\ell,M_\ell)}$ for $p\in\Pi_I$.\label{alg:mlfftinv_Il_line}
 \vspace{0.3em}
     \STATE Compute $\hat{p}^{\,\Lambda(\boldz_\ell,M_\ell)}_\boldk:=\frac{1}{M_\ell} \sum_{j=0}^{M_\ell-1} p\left(\frac{j}{M_\ell}\boldz_\ell \bmod \boldone\right) \,\mathrm{e}^{-2\pi\mathrm{i} j\boldk\cdot\boldz_\ell /M_\ell}$ for $\boldk\in I_\ell$, using adjoint rank-1 lattice FFT, cf.~\eqref{eq:r1lfft_1} and \eqref{eq:r1lfft_2}. We have $\hat{p}^{\,\Lambda(\boldz_\ell,M_\ell)}_\boldk=\hat{p}_\boldk$ for $p\in\Pi_I$.
     \label{alg:mlfftinv_lfft_line}
 \vspace{0.3em}
     \STATE For $\boldk\in I_\ell$, set $\hat{p}^{\,\Lambda}_\boldk:=\hat{p}^{\,\Lambda}_\boldk+\hat{p}^{\,\Lambda(\boldz_\ell,M_\ell)}_\boldk$ and $\mathtt{counter}[\boldk]:=\mathtt{counter}[\boldk]+1$.
     \ENDFOR
 \vspace{0.3em}
     \STATE Set $\hat{p}^{\,\Lambda}_\boldk:=\hat{p}^{\,\Lambda}_\boldk/\mathtt{counter}[\boldk]$ for $\boldk\in I$.
  \end{algorithmic}
 \vspace{1em}
 \scalebox{0.98}{%
  \begin{tabular}{p{2.2cm}p{1.5cm}p{9.6cm}}
    Output: & \rule[-0.7em]{0em}{1em}$\left(\hat{p}^{\,\Lambda}_\boldk\right)_{\boldk\in I}$ & reconstructed Fourier coefficients $=\left(\hat{p}_\boldk\right)_{\boldk\in I}$ for $p\in\Pi_I$ \\
 \end{tabular}}
 \hrule\vspace{.3em}
  \begin{tabular}{p{2.2cm}p{10cm}}
    Complexity: & $\mathcal{O}\big(M\log M + L\,|I|\,(d + \log |I|)\big)$, $M:=\sum_{\ell=1}^L M_\ell $\\
  \end{tabular}
\end{algorithm}

Similar to the approach for single rank-1 lattices $\Lambda(\boldz,M)$ in Section~\ref{sec:pre:single_r1l},
we apply Algorithm~\ref{algo:mlfft2:ifft_direct} to sampling values~$\left(f(\boldsymbol{\tilde{x}}_j)\right)_{\boldsymbol{\tilde{x}}_j\in\Lambda}$ of continuous periodic functions $f\colon\T^d\rightarrow\C$, where $\Lambda=\Lambda(\boldz_1,M_1,\ldots,\boldz_L,M_L)$ is a reconstructing multiple rank-1 lattice fulfilling reconstruction property~\eqref{equ:cond:reco_ml_mlfft2} for a given frequency index set~$I$.
We denote by $\hat{f}^{\,\Lambda}_\boldk$, $\boldk\in I$, the approximated Fourier coefficients returned by Algorithm~\ref{algo:mlfft2:ifft_direct}.
Correspondingly, we define a sampling operator by
\begin{align}
S_I^{\Lambda}f&:=\sum_{\zb k\in I}\hat{f}_{\zb k}^{\Lambda}\e^{2\pi\ii\zb k\cdot\circ}\label{eqn:sampling_operator_mr1l}.
\end{align}
The main focus of this paper are upper bounds on norms of sampling errors $f-S_I^{\Lambda}f$ as discussed in the next section,
which immediately result in the estimates of the sampling numbers~$g_M^{\mathrm{mr1l}}$ in Tables~\ref{tab:mainres:overview_general_r1l} and~\ref{tab:mainres:overview_general_r1l_2}.

\section{Multiple rank-1 lattice sampling and error estimates}
\label{sec:mult_r1l_sampl}

In this section, we show estimates on sampling errors $f - S^{\Lambda}_I f$ in specific norms,
when sampling along multiple rank-1 lattices~$\Lambda$ fulfilling reconstruction property~\eqref{equ:cond:reco_ml_mlfft2}
and applying Algorithm~\ref{algo:mlfft2:ifft_direct}.
We use the splitting approach~\eqref{equ:periodic:approx:triangle}, which leads to already estimated bounds on truncation errors $f - S_I f$, see e.g.\ \cite{ByKaUlVo16,volkmerdiss} and the references therein.
Subsequently, we determine the corresponding aliasing errors
$S_I f - S^{\Lambda}_I f$.

To this end, we introduce the index set
 $$\mathcal{L}^{\Lambda}_\boldk:=\{\ell\in\{1,\ldots,L\}\colon \boldk\cdot\boldz_\ell\not\equiv \boldh\cdot\boldz_\ell\imod{M_\ell} \; \forall \boldh\in I\setminus\{\boldk\} \}$$
for each frequency $\boldk\in I$, where the cardinality $|\mathcal{L}^{\Lambda}_\boldk|$ corresponds to
the value of $\mathtt{counter}[\boldk]$ at the end of Algorithm~\ref{algo:mlfft2:ifft_direct}.
 
\begin{Lemma}\label{lemma:aliasing_error_mlfft:exact}
 Let a function $f\in\mathcal{A}(\T^d)\cap\mathcal{C}(\T^d)$,
 a frequency index set $I\subset\Z^d$, $|I|<\infty$,
 and a multiple rank-1 lattice $\Lambda:=\Lambda(\boldz_1,M_1,\ldots,\boldz_L,M_L)\subset\T^d$
 fulfilling reconstruction property~\eqref{equ:cond:reco_ml_mlfft2}
 be given.
Then, the aliasing error $S_I f - S^{\Lambda}_I f$ can be characterized by
 \begin{equation} \label{equ:aliasing_error_mlfft:exact}
  S_I f - S^{\Lambda}_I f
  =
  - \sum_{\boldk\in I} \; \frac{1}{|\mathcal{L}^{\Lambda}_\boldk|} \sum_{\ell\in\mathcal{L}^{\Lambda}_\boldk} \sum_{\boldh\in\Lambda(\boldz_\ell,M_\ell)^\perp\setminus\{\boldzero\}} \hat{f}_{\boldk+\boldh} \; \mathrm{e}^{2\pi\mathrm{i}\boldk\cdot\circ}
 \end{equation}
 and can be estimated by
 \begin{equation}\label{equ:aliasing_error_mlfft:Linf:A}
 \Vert S_I f - S^{\Lambda}_I f | L_\infty(\T^d) \Vert \leq L\;\sum_{\boldk\in\Z^d\setminus I} |\hat{f}_{\boldk}| = L \; \Vert f - S_I f\vert \mathcal{A}(\T^d) \Vert.
 \end{equation}
\end{Lemma}

\begin{proof}
Since we have $f(\frac{j}{M_\ell}\boldz_\ell\bmod{\boldsymbol{1}}) = \sum_{\boldh\in\Z^d} \hat{f}_\boldh \, \mathrm{e}^{2\pi\mathrm{i}j\boldh\cdot\boldz_\ell/M_\ell}$, $\ell\in\{1,\ldots,L\}$, we obtain
\begin{equation*}
 \hat{f}^{\,\Lambda(\boldz_\ell,M_\ell)}_\boldk
 \overset{\eqref{equ:fhattilde}}{=}
  \sum_{\boldh\in\Lambda(\boldz_\ell,M_\ell)^\perp} \hat{f}_{\boldk+\boldh}.
\end{equation*}
Due to the reconstruction property~\eqref{equ:cond:reco_ml_mlfft2}, we obtain $\mathcal{L}^{\Lambda}_\boldk\neq \emptyset$.
Taking the averaging $\hat{f}^{\,\Lambda}_\boldk := \frac{1}{|\mathcal{L}^{\Lambda}_\boldk|} \sum_{\ell\in\mathcal{L}^{\Lambda}_\boldk} \hat{f}^{\,\Lambda(\boldz_\ell,M_\ell)}_\boldk$ in Algorithm~\ref{algo:mlfft2:ifft_direct} and the aliasing formula
\begin{align*}
  \hat{f}^{\,\Lambda}_\boldk
  &=
  \frac{1}{|\mathcal{L}^{\Lambda}_\boldk|} \sum_{\ell\in\mathcal{L}^{\Lambda}_\boldk} \sum_{\substack{\boldh\in\Z^d \\ \boldh\cdot\boldz_\ell\,\equiv\, 0 \imod{M_\ell}}} \hat{f}_{\boldk+\boldh}
  =
  \frac{1}{|\mathcal{L}^{\Lambda}_\boldk|} \sum_{\ell\in\mathcal{L}^{\Lambda}_\boldk} \sum_{\boldh\in\Lambda(\boldz_\ell,M_\ell)^\perp} \hat{f}_{\boldk+\boldh}\\
  &=\hat{f}_{\boldk}+\frac{1}{|\mathcal{L}^{\Lambda}_\boldk|} \sum_{\ell\in\mathcal{L}^{\Lambda}_\boldk} \sum_{\boldh\in\Lambda(\boldz_\ell,M_\ell)^\perp\setminus\{\boldzero\}} \hat{f}_{\boldk+\boldh}.
 \end{align*}
into account,
representation~\eqref{equ:aliasing_error_mlfft:exact} follows.
 Exploiting the properties of the index sets $I_\ell:=\{\zb k\in I\colon \nexists \zb h\in I\setminus\{\zb k\}\text{ with }\zb k\cdot \zb z_\ell\equiv\zb h\cdot \zb z_\ell\imod{M_\ell}\}$, $\ell=1,\ldots,L$, and $\mathcal{L}^\Lambda_\boldk$ yields
 \begin{align*}
 \Vert S_I f & - S^{\Lambda}_I f | L_\infty(\T^d) \Vert
 \leq
 \sum_{\boldk\in I} \; \frac{1}{|\mathcal{L}^{\Lambda}_\boldk|} \sum_{\ell\in\mathcal{L}^{\Lambda}_\boldk} \sum_{\boldh\in\Lambda(\boldz_\ell,M_\ell)^\perp\setminus\{\boldzero\}} |\hat{f}_{\boldk+\boldh}| \\
 &=
 \sum_{\ell=1}^L \sum_{\boldk\in I_\ell} \frac{1}{|\mathcal{L}^{\Lambda}_\boldk|} \sum_{\boldh\in\Lambda(\boldz_\ell,M_\ell)^\perp\setminus\{\boldzero\}} |\hat{f}_{\boldk+\boldh}|
 \le\sum_{\ell=1}^L \sum_{\boldk\in I_\ell}\sum_{\boldh\in\Lambda(\boldz_\ell,M_\ell)^\perp\setminus\{\boldzero\}} |\hat{f}_{\boldk+\boldh}|.
 \end{align*}
 It remains to show
 \begin{equation}
 \sum_{\boldk\in I_\ell} \sum_{\substack{\boldh\in\Z^d\setminus\{\boldzero\} \\ \boldh\cdot\boldz_\ell\,\equiv\, 0 \imod{M_\ell}}} |\hat{f}_{\boldk+\boldh}|=
  \sum_{\boldk\in I_\ell} \sum_{\boldh\in\Lambda(\boldz_\ell,M_\ell)^\perp\setminus\{\boldzero\}} |\hat{f}_{\boldk+\boldh}|
 \leq
 \sum_{\boldk\in\Z^d\setminus I} |\hat{f}_{\boldk}|\label{equ:aliasing_error_mlfft:ineq:trunc_error}
 \end{equation}
 for each $\ell\in\{1,\ldots,L\}$ in order to achieve the statement.
 In doing so, we fix $\ell\in\{1,\ldots,L\}$ and we now deduce
 \begin{equation}\label{equ:aliasing_error_mlfft:exact:each_freq_once}
 \boldk+\boldh\neq\boldsymbol{k'}+\boldsymbol{h'} \text{ for all } \boldk\in I_\ell, \boldsymbol{k'}\in I\setminus\{\boldk\} \text{ and } \boldh,\boldsymbol{h'}\in\Lambda(\boldz_\ell,M_\ell)^\perp
 \end{equation}
 by contradiction, similarly as in the proof of \cite[Lemma~6]{ByKaUlVo16}.
 Having the reconstruction property~\eqref{equ:cond:reco_ml_mlfft2} fulfilled,
 we assume that there exist frequencies
 $\boldk\in I_\ell$, $\boldsymbol{k'}\in I\setminus\{\boldk\}$ and $\boldh,\boldsymbol{h'}\in\Lambda(\boldz_\ell,M_\ell)^\perp$
 such that $\boldk+\boldh=\boldsymbol{k'}+\boldsymbol{h'}$.
 Then, we have $\boldk-\boldsymbol{k'}=\boldsymbol{h'}-\boldh$
 and consequently 
 $$(\boldk-\boldsymbol{k'})\cdot\boldz_\ell \equiv (\boldsymbol{h'}-\boldh)\cdot\boldz_\ell\equiv 0 \imod{M_\ell},$$
 since $\boldh,\boldsymbol{h'}\in\Lambda(\boldz_\ell,M_\ell)^\perp$ implies $\boldsymbol{h'}-\boldh\in\Lambda(\boldz_\ell,M_\ell)^\perp$, cf.~\eqref{equ:def_dual_lattice}.
 Accordingly, we obtain $(\boldk-\boldsymbol{k'})\cdot\boldz_\ell \equiv 0 \imod{M_\ell}$ or equivalently $\boldk\cdot\boldz_\ell \equiv \boldsymbol{k'}\cdot\boldz_\ell \imod{M_\ell}$
 which is in contradiction with the reconstruction property~\eqref{equ:cond:reco_ml_mlfft2}. Consequently, \eqref{equ:aliasing_error_mlfft:exact:each_freq_once} follows.
 \newline
 In particular, setting $\boldsymbol{h'}:=\boldzero$, results in $\boldk+\boldh\not\in I$ for all $\boldk\in I_\ell$ and $\boldh\in\Lambda(\boldz_\ell,M_\ell)^\perp\setminus\{\boldzero\}$.
 Moreover, due to the reconstruction property~\eqref{equ:cond:reco_ml_mlfft2} and the inclusion $I_\ell\subset I$ we have $\boldk\cdot\boldz_\ell \not\equiv \boldsymbol{k'}\cdot\boldz_\ell \imod{M_\ell}$ for all $\boldk,\boldsymbol{k'}\in I_\ell$, $\boldsymbol{k'}\neq \boldk$. Consequently, we observe that the sets 
 \begin{align*}
&\{\boldl\colon \boldl\in\Z^d\setminus\{\boldk\},\, \boldl\cdot\boldz_\ell\equiv \boldk\cdot\boldz_\ell \imod{M_\ell}\}=\\
&\quad\{\boldk+\boldh\colon \boldh\in\Lambda(\boldz_\ell,M_\ell)^\perp\setminus\{\boldzero\}\},\qquad \boldk\in I_\ell,
 \end{align*}
 are pairwise disjoint and do not contain elements from $I$, i.e.,
 $$
\bigcup_{\boldk\in I_\ell}\{\boldk+\boldh\colon \boldh\in\Lambda(\boldz_\ell,M_\ell)^\perp\setminus\{\boldzero\}\}\subset \Z^d\setminus I.
 $$
This yields \eqref{equ:aliasing_error_mlfft:ineq:trunc_error} for each $\ell\in\{1,\ldots,L\}$ and the upper bound in~\eqref{equ:aliasing_error_mlfft:Linf:A} follows.
\end{proof}

These results allow for first estimates on the sampling errors for multiple rank-1 lattice sampling.
We start with the target spaces $L_{\infty}(\T^d)$ and source spaces $\mathcal{A}^{\alpha,\beta}(\T^d)$
where the latter are subspaces of the Wiener Algebra.

\begin{Theorem} \label{theorem:sampl_error_mr1l:Linf:Aabg}
 Let a function $f\in\mathcal{A}^{\alpha,\beta}(\T^d)\cap\mathcal{C}(\T^d)$,
 the frequency index set~$I_N^{d,T}$
 and a multiple rank-1 lattice $\Lambda:=\Lambda(\boldz_1,M_1,\ldots,\boldz_L,M_L)\subset\T^d$
  fulfilling reconstruction property~\eqref{equ:cond:reco_ml_mlfft2} for $I=I_N^{d,T}$ be given,
 where the refinement $N\geq 1$, the dominating mixed smoothness $\beta\geq 0$, the isotropic smoothness $\alpha>-\beta$, and the shape parameter $T:=-\frac{\alpha}{\beta}$ with $T:=-\infty$ for $\beta=0$. %
 Moreover, let the approximated Fourier coefficients~$\hat{f}^{\,\Lambda}_\boldk$, $\boldk\in I_N^{d,T}$, be computed by Algorithm~\ref{algo:mlfft2:ifft_direct}.
 Then, the sampling error is bounded by
 \begin{equation} %
  \Vert f - S^{\Lambda}_{ I_N^{d,T}} f\vert L_\infty(\T^d) \Vert %
  \leq
  \Vert f - S^{\Lambda}_{ I_N^{d,T}} f\vert \mathcal{A}(\T^d) \Vert
  \leq
  N^{-(\alpha+\beta)} \, (1+ L) \, \Vert f\vert \mathcal{A}^{\alpha,\beta}(\T^d) \Vert.
  \label{equ:sampl_error_mr1l:Linf:Aabg}
 \end{equation}
\end{Theorem}
\begin{proof}
Applying inequality~\eqref{equ:periodic:approx:triangle} in the $L_\infty(\T^d)$ norm on $f - S^{\Lambda}_{ I_N^{d,T}} f$, we estimate the sampling error by
 $\Vert f - S^{\Lambda}_{ I_N^{d,T}} f\vert L_\infty(\T^d) \Vert \leq \Vert f - S_{ I_N^{d,T}} f\vert L_\infty(\T^d) \Vert + \Vert S_{ I_N^{d,T}} f - S^{\Lambda}_{ I_N^{d,T}} f\vert L_\infty(\T^d) \Vert.$
\newline
In the following, we estimate the truncation error $f - S_{ I_N^{d,T}} f$ as in the proof of \cite[Theorem~3.3]{KaPoVo13}.
We have $f - S_{ I_N^{d,T}} f = \sum_{\boldk\in \Z^d\setminus I_N^{d,T}} \hat{f}_\boldk \,  \mathrm{e}^{2\pi\mathrm{i}\boldk\circ}$.
Using the weights~$\omega^{\alpha,\beta}(\boldk)$,
the definition of the frequency index sets~$I_N^{d,T}$ in~\eqref{equ:I_N_d_T}, and the choice $T=-\alpha/\beta$ for $\beta>0$,
we obtain
 \begin{align*}
 \Vert f - S_{ I_N^{d,T}} f\vert L_\infty(\T^d) \Vert
 &\leq
 \Vert f - S_{ I_N^{d,T}} f\vert \mathcal{A}(\T^d)\Vert
 =
 \sum_{\boldk\in \Z^d\setminus I_N^{d,T}} \vert\hat{f}_\boldk\vert \\
& =
 \sum_{\boldk\in \Z^d\setminus I_N^{d,T}} \frac{\omega^{\alpha,\beta}(\boldk)}{\omega^{\alpha,\beta}(\boldk)} \; \vert\hat{f}_\boldk\vert \\
 &\leq
\sup_{\boldk\in \Z^d\setminus I_N^{d,T}} \frac{1}{\omega^{\alpha,\beta}(\boldk)} \;\sum_{\boldk\in \Z^d\setminus I_N^{d,T}} \omega^{\alpha,\beta}(\boldk) \; \vert\hat{f}_\boldk\vert \\
& \leq
\frac{1}{\sup_{\boldk \in I_N^{d,T}}\omega^{\alpha,\beta}(\boldk)} \;\sum_{\boldk\in \Z^d} \omega^{\alpha,\beta}(\boldk) \; \vert\hat{f}_\boldk\vert \\
& \leq
 N^{-(\alpha+\beta)} \, \Vert f\vert \mathcal{A}^{\alpha,\beta}(\T^d) \Vert,
 \end{align*}
due to H\"older's inequality, see also \cite[Lemma~3.2]{KaPoVo13}.
 For $\beta=0$, we proceed analogously.\newline
 Next, we estimate the aliasing error.
 Applying Lemma~\ref{lemma:aliasing_error_mlfft:exact} results in the estimate
 $
   \Vert S_{ I_N^{d,T}} f - S^{\Lambda}_{ I_N^{d,T}} f\vert L_\infty(\T^d) \Vert
   \leq
   L \,
   \sum_{\boldk\in \Z^d\setminus I_N^{d,T}} \vert\hat{f}_{\boldk}\vert
 $
 and we proceed as in the estimate of the truncation error $f - S_{ I_N^{d,T}} f$, %
 which yields~\eqref{equ:sampl_error_mr1l:Linf:Aabg}.
\end{proof}

The last theorem provides estimates on the sampling error in terms of the refinement~$N$ of the used
frequency index sets~$I_N^{d,T}$ and the number~$L$ of used rank-1 lattices.
Since the main objectives of the paper are estimates on sampling errors
in terms of the number~$|\Lambda|$ of used sampling values, we still need to bound
$N$ and $L$ with respect to $M:=\sum_{\ell=1}^{L}M_\ell\ge|\Lambda|$.

Due to \cite[Corollary~3.7]{Kae17} a multiple rank-1 lattice~$\Lambda$
fulfilling reconstruction proper\-ty~\eqref{equ:cond:reco_ml_mlfft2} can be determined for arbitrary index sets $I\subset\Z^d$
with $|\Lambda|\lesssim  |I|\,L\lesssim |I|\,\ln|I|$, where one assumes
$$|I|\ge2,\quad
|I|\gtrsim N_I:=\max_{j=1,\ldots, d}\{\max_{\zb k\in I}k_j-\min_{\zb l\in I}l_j\},\quad\text{and}\quad
|I|\gtrsim L\,\ln L.$$
Algorithm~\ref{alg:construct_mr1l_I_distinct_primes}
efficiently constructs such a multiple rank-1 lattice~$\Lambda$
with high probability.
We stress the fact, that for $I:=I_N^{d,T}$ the first two assumptions
$|I_N^{d,T}|\ge2$ and $|I_N^{d,T}|\ge N_I$ are naturally fulfilled.
The third assumption $|I_N^{d,T}|\gtrsim L\,\ln L$ is not a restriction in asymptotics
since $L$ depends only linearly on $\ln|I_N^{d,T}|$.

Accordingly, we achieve the following statements that lead
to the results of the second columns and third rows of Tables \ref{tab:mainres:overview_general_r1l} and
\ref{tab:mainres:overview_general_r1l_2} when taking Lemma \ref{lemma:estimate:Aabg:Habg} into account.

\begin{Corollary}\label{cor:sampl_error_mr1l:Linf:Aabg}
Under the assumptions of Theorem~\ref{theorem:sampl_error_mr1l:Linf:Aabg}, we additionally assume
that the number $L$ of single rank-1 lattices is bounded by $L\lesssim\log |I_N^{d,T}|$ and that the number of samples~$M$ is bounded
by $M:=\sum_{\ell=1}^{L}M_\ell \lesssim |I_N^{d,T}| \log |I_N^{d,T}|$
 (e.g.~$\Lambda$ is constructed by Algorithm~\ref{alg:construct_mr1l_I_distinct_primes}).
 Then the sampling error is bounded by
 \begin{align*}
   \Vert f & - S^{\Lambda}_{ I_N^{d,T}} f\vert L_\infty(\T^d) \Vert
   \leq
   \Vert f - S^{\Lambda}_{ I_N^{d,T}} f\vert \mathcal{A}(\T^d) \Vert
   \\
   &\lesssim
   \Vert f\vert \mathcal{A}^{\alpha,\beta}(\T^d) \Vert \;
   (\log M)
   \begin{cases}
   \left(\frac{\log M}{M}\right)^{\frac{\alpha}{d}+\beta} & \text{for } \alpha>0,\,\beta\geq 0,\\
   \left(\frac{\log^d M}{M}\right)^{\beta} & \text{for } \beta>\alpha=0,\\
   \left(\frac{\log M}{M}\right)^{\alpha+\beta} & \text{for } -\beta<\alpha<0,
   \end{cases}
 \end{align*}
 where the constants may depend on the dimension~$d$ and the smoothness parameters $\alpha,\beta$.
\end{Corollary}
\begin{proof}
The assumed upper bounds on $L$ and $M$ imply $L\lesssim \log N\lesssim \log M$, cf.\ Lemma~\ref{lemma:index_set:omega}. 
Theorem~\ref{theorem:sampl_error_mr1l:Linf:Aabg} yields
 \begin{equation*}
 \Vert f - S^{\Lambda}_{ I_N^{d,T}} f\vert L_\infty(\T^d) \Vert
 \leq
 \Vert f - S^{\Lambda}_{ I_N^{d,T}} f\vert \mathcal{A}(\T^d) \Vert
 \lesssim
 N^{-(\alpha+\beta)} \, (\log N) \, 
\Vert f\vert \mathcal{A}^{\alpha,\beta}(\T^d) \Vert,
 \end{equation*}
 where the constants may depend on the dimension~$d$ and the shape parameter~$T$.
Again, taking the cardinalities $|I_N^{d,T}|$ from Lemma~\ref{lemma:index_set:omega} into account yields the assertion.
\end{proof}

Since users are  also interested in approximations in more specific target spaces, cf.\ e.g.\ \cite{Ys10, GriHa13},
we investigate approximation errors measured in the norm of the spaces $\mathcal{H}^{r,t}(\T^d)$.
In this setting, commonly occurring source spaces are Hilbert spaces of the same type with smoothness of higher order.
Estimates on the norm of the truncation error $f - S_{ I_N^{d,T}}$ can be obtained by an already established argument which uses shiftings of the smoothness parameters.
However, this shifting argument is not applicable in order to achieve suitable estimates on the aliasing errors $S_{ I_N^{d,T}} f - S^{\Lambda}_{ I_N^{d,T}}$ of rank-1 lattice sampling strategies.
For that reason, we take a detour through subspaces of the Wiener Algebra in order to estimate the norm of the aliasing error, which leads to the
crucial statement in the next theorem.

\begin{Theorem} \label{theorem:approx_error:Hrtg:Aabg}
 Let the dominating mixed smoothness $t\geq 0$ and the isotropic smoothness $r\geq -t$ of the target space~$\mathcal{H}^{r,t}(\T^d)$
 as well as the dominating mixed smoothness $\beta\geq t\geq 0$ and the isotropic smoothness $\alpha>-\beta$ of the source space~$\mathcal{A}^{\alpha,\beta}(\T^d)$ be given,
 where $\alpha+\beta > r+t\ge 0$.
 Moreover, let
 a function $f\in\mathcal{C}(\T^d)\cap\mathcal{A}^{\alpha,\beta}(\T^d)$,
 a frequency index set $ I_N^{d,T}$
 and a multiple rank-1 lattice $\Lambda:=\Lambda(\boldz_1,M_1,\ldots,\boldz_L,M_L)\subset\T^d$
 fulfilling reconstruction property~\eqref{equ:cond:reco_ml_mlfft2} for $I=I_N^{d,T}$ be given,
 where the refinement $N\geq 1$
 and the shape parameter
 $T\in [-\frac{r}{t},1)$ with $-\frac{r}{t}:=-\infty$ for $t=0$.
 Then, the sampling error is bounded by
 \begin{align} \nonumber
  \Vert f -  S^{\Lambda}_{ I_N^{d,T}} f\vert & \mathcal{H}^{r,t}(\T^d) \Vert
  \leq
  N^{-(\alpha-r+\beta - t)}
  \\
  \cdot \nonumber
  & \left(
  \Vert f\vert \mathcal{H}^{\alpha,\beta}(\T^d) \Vert
  \left\{
  \begin{array}{ll}
   N^{(d-1)\frac{T(\beta - t)+\alpha-r}{d-T}} & \textrm{ for } T > -\frac{\alpha-r}{\beta-t} \\
   d^{-\frac{T(\beta - t) + \alpha-r}{1-T}} & \textrm{ for } T \leq -\frac{\alpha-r}{\beta-t}
  \end{array}
  \right\}
  \right. \\
  &\left. +
     \Vert f\vert \mathcal{A}^{\alpha,\beta}(\T^d) \Vert \, L \,
  \left\{
  \begin{array}{ll}
      d^{\frac{Tt+r}{1-T}} \, N^{(d-1)\frac{T\beta+\alpha}{d-T}} & \textrm{ for } T > -\frac{\alpha}{\beta} \\
      d^{-\frac{T(\beta-t) + \alpha-r}{1-T}} & \textrm{ for } T \leq -\frac{\alpha}{\beta}
  \end{array}
  \right\}
  \right), \label{equ:approx_error_mlfft:Hrtg:Aabg}
 \end{align}
   where we define $T t:=0$ for $t=0$. 
\end{Theorem}
\begin{proof}
We follow the general strategy of the proof of \cite[Theorem~3.4]{KaPoVo13}.
First, we apply inequality~\eqref{equ:periodic:approx:triangle} in the $\mathcal{H}^{r,t}(\T^d)$ norm on $f - S^{\Lambda}_{ I_N^{d,T}} f$
and split up the sampling error.
For a function $f\in\mathcal{H}^{\alpha,\beta}(\T^d)$,
we have
  \begin{eqnarray*}
  \Vert f &-& S_{ I_N^{d,T}} f\vert \mathcal{H}^{r,t}(\T^d) \Vert^2
  =
  \sum_{\boldk\in \Z^d\setminus I_N^{d,T}} \omega^{r,t}(\boldk)^2 \; \frac{\omega^{\alpha,\beta}(\boldk)^2}{\omega^{\alpha,\beta}(\boldk)^2} \; \vert\hat{f}_\boldk\vert^2 \\
  &\leq&
  \left(
  \max_{\boldk\in \Z^d\setminus I_N^{d,T}}
  \omega^{-(\alpha-r),-(\beta-t)}(\boldk)^2
  \right)
  \sum_{\boldk\in \Z^d\setminus I_N^{d,T}}
  \omega^{\alpha,\beta}(\boldk)^2 \; \vert\hat{f}_\boldk\vert^2.
  \end{eqnarray*}
  Applying \cite[Lemma~3.2]{KaPoVo13} with $\tilde{\alpha}:=\alpha-r$ and $\tilde{\beta}:=\beta-t$ yields
  $$
   \max_{\boldk\in \Z^d\setminus I_N^{d,T}} \omega^{-(\alpha-r),-(\beta-t)}(\boldk) \leq N^{-(\alpha-r+\beta - t)} 
     \begin{cases}
      N^{(d-1)\frac{T(\beta - t)+\alpha-r}{d-T}} & \textrm{for } T > -\frac{\alpha-r}{\beta-t}, \\
      d^{-\frac{T(\beta - t) + \alpha-r}{1-T}} & \textrm{for } T \leq -\frac{\alpha-r}{\beta-t},
     \end{cases}
  $$
  and consequently the first summand of the right hand side in~\eqref{equ:approx_error_mlfft:Hrtg:Aabg} follows.\newline
  Next, we estimate the aliasing error.
  Since we have the aliasing formula~\eqref{equ:aliasing_error_mlfft:exact} in Lemma~\ref{lemma:aliasing_error_mlfft:exact}
  and due to the concaveness of the square root function, we estimate
   \begin{align}
    \Vert S_{ I_N^{d,T}} f - S^{\Lambda}_{ I_N^{d,T}} f\vert \mathcal{H}^{r,t}(\T^d) \Vert
    &\leq \nonumber
    \left(
    \sum_{\boldk\in I_N^{d,T}} \frac{\omega^{r,t}(\boldk)^2}{|\mathcal{L}^{\Lambda}_\boldk|^2} \left\vert \sum_{\ell\in\mathcal{L}^{\Lambda}_\boldk} \sum_{
    \boldh\in\Lambda(\boldz_\ell,M_\ell)^\perp\setminus\{\boldzero\}
    } \hat{f}_{\boldk+\boldh} \right\vert^2
    \right)^{\frac{1}{2}} \\
    &\leq \nonumber
    \sum_{\boldk\in I_N^{d,T}} \frac{\omega^{r,t}(\boldk)}{|\mathcal{L}^{\Lambda}_\boldk|} \left\vert \sum_{\ell\in\mathcal{L}^{\Lambda}_\boldk} \sum_{
    \boldh\in\Lambda(\boldz_\ell,M_\ell)^\perp\setminus\{\boldzero\}
    } \hat{f}_{\boldk+\boldh} \right\vert \\
    &\leq \nonumber
    \max_{\boldk\in I_N^{d,T}} \omega^{r,t}(\boldk) \, 
    \sum_{\boldk\in I_N^{d,T}} \sum_{\ell\in\mathcal{L}^{\Lambda}_\boldk} \sum_{
    \boldh\in\Lambda(\boldz_\ell,M_\ell)^\perp\setminus\{\boldzero\}
    } \left\vert \hat{f}_{\boldk+\boldh} \right\vert \\
    &=\nonumber
    \max_{\boldk\in I_N^{d,T}} \omega^{r,t}(\boldk) \, 
    \sum_{\ell=1}^L \sum_{\boldk\in I_\ell} \sum_{\boldh\in\Lambda(\boldz_\ell,M_\ell)^\perp\setminus\{\boldzero\}} |\hat{f}_{\boldk+\boldh}| \\
    &\overset{\eqref{equ:aliasing_error_mlfft:ineq:trunc_error}}{\leq}\nonumber
    \max_{\boldk\in I_N^{d,T}} \omega^{r,t}(\boldk) \;
     L\, \sum_{\boldk\in\Z^d\setminus I_N^{d,T}} |\hat{f}_{\boldk}|.
   \end{align}
  Applying \cite[Lemma~2.4]{KaPoVo13} with $-\frac{r}{t}=:\tilde{T}\le T$, we estimate
   $\max_{\boldk\in I_{N}^{d,T}}\omega^{r,t}(\boldk)\leq d^{(Tt+r)/(1-T)}N^{r+t}$,
   where we define $T t:=0$ for $t=0$.
  Incorporating the weights $\omega^{\alpha,\beta}(\boldk)$, we infer
\begin{align*}
 \Vert S_{ I_N^{d,T}} f - S^{\Lambda}_{ I_N^{d,T}} f\vert \mathcal{H}^{r,t}(\T^d) \Vert
 &\leq
 d^{\frac{Tt+r}{1-T}}N^{r+t} L\, \max_{\boldk\in\Z^d\setminus I_N^{d,T}} \frac{1}{\omega^{\alpha,\beta}(\boldk)} \, \Vert f\vert \mathcal{A}^{\alpha,\beta}(\T^d) \Vert. %
\end{align*}
Again we apply \cite[Lemma~3.2]{KaPoVo13}, where we use $\tilde{\alpha}:=\alpha$ and $\tilde{\beta}:=\beta$.
Thus, the latter term can be bounded from above by
\begin{align*}
 \Vert f\vert \mathcal{A}^{\alpha,\beta}(\T^d) \Vert\,
  N^{-(\alpha-r+\beta - t)}\, L\,
  \begin{cases}
      d^{\frac{Tt+r}{1-T}} \, N^{(d-1)\frac{T\beta+\alpha}{d-T}} & \textrm{ for } T > -\frac{\alpha}{\beta}, \\
      d^{-\frac{T(\beta-t) + \alpha-r}{1-T}} & \textrm{ for } T \leq -\frac{\alpha}{\beta},
  \end{cases}
\end{align*}
and we obtain~\eqref{equ:approx_error_mlfft:Hrtg:Aabg}.
\end{proof}

Note that Theorem \ref{theorem:sampl_error_mr1l:Linf:Aabg} as well as Theorem \ref{theorem:approx_error:Hrtg:Aabg} provide estimates on the
sampling errors with respect to the refinement $N$ of the frequency index sets~$I_N^{d,T}$, incorporating the dependencies on all parameters
and keeping track of all constants.
For determining sampling numbers, we consider the dependency on the number $M$  of used sampling values without regarding terms that only
depend on the fixed smoothness parameters $r$,$t$,$\alpha$,$\beta$, the shape parameter $T$, or the fixed dimension $d$.
To this end, we apply embedding arguments for the considered function spaces.

In the following, we only discuss these parameter combinations that lead to best possible main rates of the approximation errors with respect to the refinement $N$ in~\eqref{equ:approx_error_mlfft:Hrtg:Aabg}. 

\begin{Corollary}\label{cor:approx_error:Hrtg:Aabg}
Under the assumptions of Theorem~\ref{theorem:approx_error:Hrtg:Aabg}, we choose the shape parameter
$T\in \left[-\frac{r}{t},\min\{ -\frac{\alpha}{\beta}, -\frac{\alpha-r}{\beta-t}\}\right]$.
Then the sampling error is bounded by
 \begin{align*}
  \Vert f - S^{\Lambda}_{ I_N^{d,T}} f\vert \mathcal{H}^{r,t}(\T^d) \Vert
  &\lesssim
  N^{-(\alpha-r+\beta - t)} \, L \, \Vert f\vert \mathcal{A}^{\alpha,\beta}(\T^d) \Vert,
 \end{align*}
 where the (hidden) constants may depend on the parameters $d,r,t,\alpha,\beta,T$.
\end{Corollary}

Using the embeddings $\mathcal{H}^{\alpha,\beta}(\T^d) \hookrightarrow \mathcal{A}^{\alpha,\beta-\lambda}(\T^d)$, $\lambda>1/2$, from Lemma~\ref{lemma:estimate:Aabg:Habg} immediately yields
$\mathcal{H}^{r,t}(\T^d)$-error estimates for continuous functions belonging to Sobolev spaces~$\mathcal{H}^{\alpha,\beta}(\T^d)$ of generalized mixed smoothness.
Furthermore, we obtain sampling numbers for multiple rank-1 lattice sampling using the estimates on the cardinality of the frequency index sets~$I_N^{d,T}$, which are stated in Lemma \ref{lemma:index_set:omega}, and the estimate $L\lesssim \log|I_N^{d,T}|\lesssim \log M$.

\begin{Corollary}\label{cor:approx_error:Hrtg:Habg}
 Let smoothness parameters $\beta\ge\lambda>1/2$, $r,t,\alpha\in\R$, $\beta-\lambda\geq t\geq 0$, $\alpha+\beta-\lambda > r+t\ge 0$,
 the shape parameter $T\in \left[-\frac{r}{t},\min\{ -\frac{\alpha}{\beta-\lambda}, -\frac{\alpha-r}{\beta-\lambda-t}\}\right]$
  with $-\frac{r}{t}:=-\infty$ for $t=0$,
 a function $f\in\mathcal{H}^{\alpha,\beta}(\T^d)\cap\mathcal{C}(\T^d)$,
 a frequency index set $ I_N^{d,T}$
 and a multiple rank-1 lattice $\Lambda:=\Lambda(\boldz_1,M_1,\ldots,\boldz_L,M_L)\subset\T^d$
 fulfilling reconstruction property~\eqref{equ:cond:reco_ml_mlfft2} for $I=I_N^{d,T}$ be given,
 where the refinement $N\geq 1$.
Additionally, let $L\lesssim\log |I_N^{d,T}|$ and $M:=\sum_{\ell=1}^{L}M_\ell \lesssim |I_N^{d,T}| \log |I_N^{d,T}|$,
(e.g.~$\Lambda$ is constructed by Algorithm~\ref{alg:construct_mr1l_I_distinct_primes}).
Then the sampling error is bounded by
 \begin{align*}
  \Vert f - S^{\Lambda}_{ I_N^{d,T}} f\vert \mathcal{H}^{r,t}(\T^d) \Vert
  &\lesssim
  N^{-(\alpha-r+\beta - t - \lambda)} \, (\log N) \, \Vert f\vert \mathcal{H}^{\alpha,\beta}(\T^d) \Vert
 \end{align*}
 and
 \begin{align*}
  \Vert f -& S^{\Lambda}_{ I_N^{d,T}} f\vert \mathcal{H}^{r,t}(\T^d) \Vert
  \\&\lesssim
  \Vert f\vert \mathcal{H}^{\alpha,\beta}(\T^d) \Vert
   \,
   (\log M)
   \begin{cases}
   \left(\frac{\log M}{M}\right)^{\frac{\alpha-r+\beta-\lambda}{d}} & \text{for } T=-\infty,\\
   \left(\frac{\log M}{M}\right)^{\frac{T/d-1}{T-1}(\alpha-r+\beta - t -\lambda)} & \text{for } -\infty<T<0\\
   \left(\frac{\log^d M}{M}\right)^{\alpha - r + \beta - t -\lambda} & \text{for } T=0,\\
   \left(\frac{\log M}{M}\right)^{\alpha-r+\beta - t-\lambda} & \text{for } T>0,
   \end{cases}
 \end{align*}
 where the constants may depend on the dimension $d$ and the smoothness parameters $r,t,\alpha,\beta,\lambda$.
\end{Corollary}

\begin{table}[ht]
\centering
\begin{tabular}{lllll}
\toprule
\multicolumn{2}{l}{possible} & \multicolumn{2}{l}{} & largest \\
\multicolumn{2}{l}{parameter choices} & \multicolumn{2}{l}{parameter restrictions} & possible $T$\\
\midrule
$r=0$ & $\beta -\lambda \ge t = 0$ & $T \in [-\infty,-\frac{\alpha}{\beta-\lambda}]$ & $\alpha\in\R$ & $-\frac{\alpha}{\beta-\lambda}$  \\
$r>0$ & $\beta -\lambda \ge t = 0$ & $T \in [-\infty,-\frac{\alpha}{\beta-\lambda}]$ & $\alpha\in\R$ & $-\frac{\alpha}{\beta-\lambda}$  \\
$r=0$ & $\beta -\lambda \ge t > 0$ & $T \in [0,-\frac{\alpha}{\beta-\lambda}]$ & $\alpha \le 0$ & $-\frac{\alpha}{\beta-\lambda}\ge 0$  \\
$r>0$ & $\beta -\lambda \ge t > 0$ & $T \in [-\frac{r}{t},-\frac{\alpha}{\beta-\lambda}]$  & $\alpha\le\frac{r}{t}(\beta-\lambda)$ & $-\frac{\alpha}{\beta-\lambda}$  \\
$r<0$ & $\beta -\lambda \ge t > 0$ & $T \in [-\frac{r}{t},-\frac{\alpha}{\beta-\lambda}]$ & $\alpha\le\frac{r}{t}(\beta-\lambda) < 0$ & $-\frac{\alpha}{\beta-\lambda}>0$  \\
\bottomrule

\end{tabular}
\caption{Possible parameter combinations and corresponding restrictions for Corollary~\ref{cor:approx_error:Hrtg:Habg}. Note the additional assumptions
$\beta\ge\lambda>1/2$ and $\alpha+\beta-\lambda > r+t\ge 0$.}\label{table:parameter:cor:approx_error:Hrtg:Habg}
\end{table}

\begin{Remark}
The given intervals for $T$ are crucial restrictions of Corollaries~\ref{cor:approx_error:Hrtg:Aabg} and~\ref{cor:approx_error:Hrtg:Habg}.
In particular, there exist parameter combinations, where the intervals are empty and thus
the corollaries do not help. Nevertheless, many neglected parameter
combinations can be treated by applying Theorem~\ref{theorem:approx_error:Hrtg:Aabg}.
However, we only restricted the considerations to the cases where we get error rates that are best possible with respect to the refinement~$N$.
Consequently, the assertions only hold for
carefully chosen parameter combinations.
\newline
Table~\ref{table:parameter:cor:approx_error:Hrtg:Habg} presents these combinations for Corollary~\ref{cor:approx_error:Hrtg:Habg}.
For $\alpha\ge 0$, it turns out that all adequate parameter combinations imply that $T:=-\frac{\alpha}{\beta-\lambda}$
is the best possible choice of this shape parameter in the sense that for $T\in[-\infty,0]$ the dependency on the dimension $d$ of the sampling rates decreases for growing $T$. Accordingly, the frequency index set $I_N^{d,T}$
only depends on $\alpha$ and $\beta-\lambda$, and thus, the smoothness parameters $\alpha$ and $\beta$ of the source space determine the
estimate of the number of sampling nodes $M$ in terms of the refinement~$N$.
For $\alpha<0$, the shape parameter $T$ can be arbitrarily chosen in the interval $\left(0,-\frac{\alpha}{\beta-\lambda}\right]\cap\left[-\frac{r}{t},-\frac{\alpha}{\beta-\lambda}\right]$ in order to achieve the best possible statement with respect to the refinement $N$ in Corollary~\ref{cor:approx_error:Hrtg:Habg}.
\newline
Moreover, it may be helpful to increase the offset $\lambda$ in order to get approximation estimates from Corollary~\ref{cor:approx_error:Hrtg:Habg}.
For, e.g., $r=-1$, $t=2$, $\alpha=-2$, $\beta=5$, the offset $\lambda$ near 1/2 will not result in appropriate parameter constellations, since
the restriction $\alpha\le \frac{r}{t}(\beta-\lambda)$ is violated.
Increasing the offset~$\lambda$ to at least one will retrieve this condition.
In this way Corollary~\ref{cor:approx_error:Hrtg:Habg} also yields estimates for this parameter combination even
though these estimates are far away from optimal ones.
\end{Remark}

\begin{Remark}
In contrast to the approximation estimates concerning sampling based on single rank-1 lattices, cf.\ columns three in Tables~\ref{tab:mainres:overview_general_r1l} and \ref{tab:mainres:overview_general_r1l_2},
higher dominating mixed smoothness $\beta$ of the function under consideration results in
upper bounds of the sampling error that behave best possible with respect to the increase of $\beta$.
However, comparing to the known lower and the best possible upper bounds, 
the sampling rates $g_M^{\mathrm{mr1l}}$ for multiple rank-1 lattices suffer from an offset $\lambda$ in the exponent of the main terms that is slightly larger than 1/2, cf.\ columns four in Tables~\ref{tab:mainres:overview_general_r1l} and \ref{tab:mainres:overview_general_r1l_2}.
This offset is caused by the proof strategy that uses embeddings
of the Hilbert spaces of generalized mixed smoothness in subspaces of the Wiener algebra.
In the next section, we present numerical tests where we do not recognize
the offset~$\lambda$.  This is one of the reasons why we conjecture that optimal upper bounds
of the sampling errors for the source spaces $\mathcal{H}^{\alpha,\beta}(\T^d)$
\pagebreak
do not suffer from this offset.
However, further improvements on the assertions in Corollary~\ref{cor:approx_error:Hrtg:Habg} require completely different proof techniques, which are
--- at least for the authors ---
currently not available.
\end{Remark}

\section{Numerical results}\label{sec:numerics}

In this section, we perform numerical tests in up to 20 dimensions. 

\subsection[Test functions using symmetric hyperbolic cross frequency index sets]{Test functions $G_{3,4}^d$ and $G_3^d$ using symmetric hyperbolic cross frequency index sets $I_N^{d,0}$}

We consider the tensor-product test functions
$G_{3,4}^d\colon\T^d\rightarrow\C$ from~\cite{KaPoVo13} and~\cite{KaPoVo14},
$G_{3,4}^d(\boldx):=\prod_{s=1}^{d} g_{3,4}(x_s)$,
where the one-dimensional function $g_{3,4}\colon\T\rightarrow\C$ is defined by
$$
g_{3,4}(x):=C \left(4+\sgn((x\bmod 1)-1/2) \left[\sin(2\pi x)^3 + \sin(2\pi x)^4\right]\right).
$$
The constant $C=8 \sqrt{\frac{6\pi}{6369\pi - 4096}}$ is fixed
and
$\sgn$ denotes the sign function, $\sgn(x):=x/\vert x\vert$ for $x\neq 0$ and $\sgn(0):=0$.
We have
$\Vert G_{3,4}^d \vert L^2(\T^d) \Vert = 1$
as well as
$G_{3,4}^d\in\mathcal{A}^{0,3-\epsilon}(\T^d)$, $G_{3,4}^d\in\mathcal{H}^{0,\frac{7}{2}-\epsilon}(\T^d)$, $\epsilon > 0$.
Moreover, we remark that $\Vert G_{3,4}^d \vert \mathcal{A}(\T^d) \Vert = \left(8 \sqrt{\frac{6 \pi }{6369 \pi -4096}} \left(4+\frac{388}{105 \pi }\right)\right)^d\approx (1.42522)^d$.

We approximate the test functions~$G_{3,4}^d$ by multivariate trigonometric polynomials $p$ based on samples along reconstructing single and reconstructing multiple rank-1 lattices, i.e., using the rank-1 lattice sampling operators~$S_I^{\Lambda(\zb z,M)}$ and~$S_I^{\Lambda}$. As frequency index sets $I$, we use symmetric hyperbolic crosses $I_N^{d,0}$ with various refinements $N\in\N$ in dimensions $d\in\{2,\ldots,8\}$.
We use reconstructing single rank-1 lattices~$\Lambda(\boldz,M)$ generated by the implementation \cite[\texttt{genlattice\_cbc\_incr\_bisect}]{Vo_taylorR1Lnfft} of \cite[Algorithm~3.7]{kaemmererdiss}, see also \cite[Table~6.2]{KaPoVo13}. The reconstructing multiple rank-1 lattices~$\Lambda$ are generated by Algorithm~\ref{alg:construct_mr1l_I_distinct_primes}.

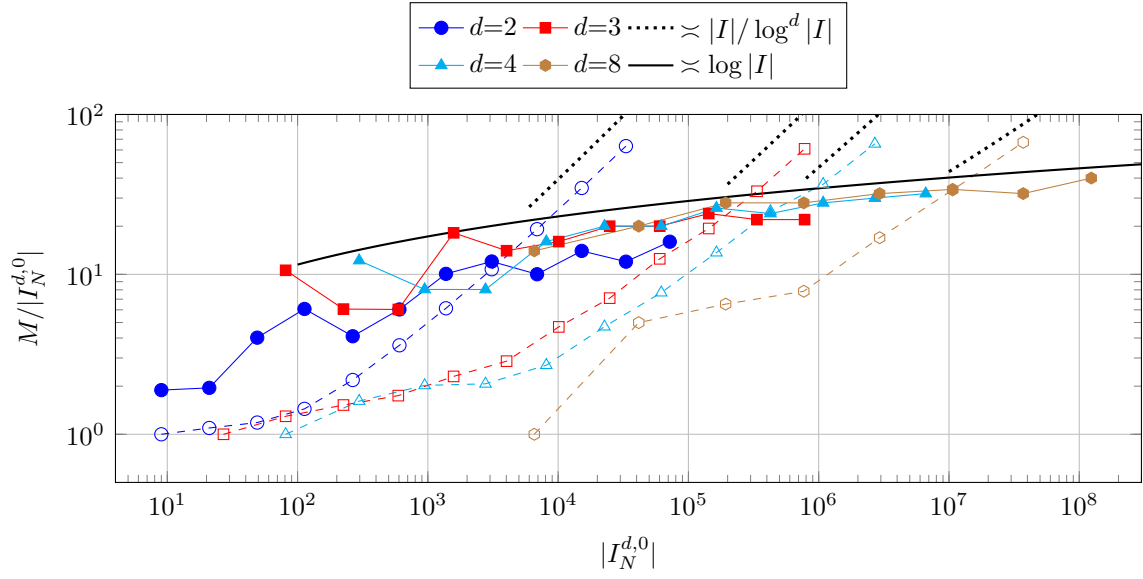
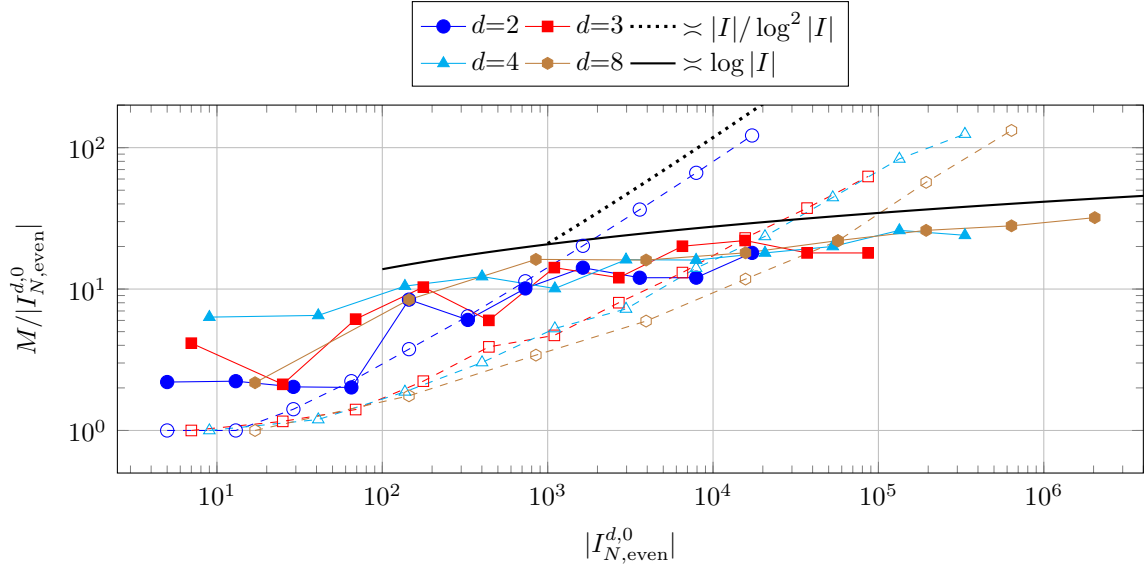
\begin{figure}[htb]
\hfill
\subfloat[$I:=I_{N}^{d,0}$]{\label{fig:numerics:periodic:reco_r1l:oversampling:hc}
\begin{tikzpicture}
  \begin{loglogaxis}[enlargelimits=false,xmin=4,xmax=3e8,ymin=0.5,ymax=1e2,ytick={1,10,100},height=0.4\textwidth, width=0.94\textwidth, grid=major, xlabel={$|I_{N}^{d,0}|$}, ylabel={$M/|I_{N}^{d,0}|$}, font=\footnotesize,
  legend style={at={(0.5,1.07)}, anchor=south,legend columns=3,legend cell align=left, font=\footnotesize, %
  },]
\addplot[forget plot,dashed,mark options={solid},blue,mark=o,mark size=2.5] coordinates {
 (9,1.000) (21,1.095) (49,1.184) (113,1.442) (265,2.185) (605,3.602) (1377,6.137) (3093,10.761) (6889,19.175) (15169,34.698) (33145,63.396)
};
\addplot[blue,mark=*,mark size=2.5] coordinates {
 (9,1.889) (21,1.952) (49,4.020) (113,6.080) (265,4.102) (605,6.035) (1377,10.075) (3093,12.041) (6889,10.016) (15169,14.019) (33145,12.007) (71941,16.008)
};
\addlegendentry{$d$=2}
\addplot[forget plot,dashed,mark options={solid},red,mark=square,mark size=2] coordinates {
 (27,1.000) (81,1.296) (225,1.524) (593,1.745) (1577,2.301) (4021,2.866) (10113,4.693) (24869,7.101) (60217,12.509) (143225,19.367) (336033,33.097) (779557,60.975)
};
\addplot[red,mark=square*,mark size=2] coordinates {
 (81,10.605) (225,6.067) (593,6.035) (1577,18.186) (4021,14.057) (10113,16.033) (24869,20.015) (60217,20.018) (143225,24.003) (336033,22.002) (779557,22.001)
};
\addlegendentry{$d$=3}
\addplot[forget plot,black,domain=6e3:1e5,samples=100,dotted,very thick] {x/(ln(x))^2/3};
\addplot[forget plot,black,domain=2e5:1e8,samples=100,dotted,very thick] {x/(ln(x))^3/3};
\addplot[forget plot,black,domain=8e5:1e8,samples=100,dotted,very thick] {x/(ln(x))^4*1.7};
\addplot[black,domain=1e7:1e8,samples=100,dotted,very thick] {x/(ln(x))^8*20000};
\addlegendentry{$\asymp |I|/\log^d|I|$}
\addplot[forget plot,dashed,mark options={solid},cyan,mark=triangle,mark size=2.5] coordinates {
 (81,1.000) (297,1.613) (945,2.022) (2769,2.068) (8113,2.705) (22665,4.708) (61889,7.688) (164137,13.672) (426193,24.781) (1082305,36.619) (2698849,65.527)
};
\addplot[cyan,mark=triangle*,mark size=2.5] coordinates {
 (297,12.205) (945,8.050) (2769,8.039) (8113,16.038) (22665,20.018) (61889,20.008) (164137,26.005) (426193,24.003) (1082305,28.002) (2698849,30.000) (6630233,32.000)
};
\addlegendentry{$d$=4}
\addplot[forget plot,dashed,mark options={solid},brown,mark=hexagon,mark size=2.25] coordinates {
 (6561,1.000) (41553,4.991) (193185,6.518) (768609,7.843) (2935521,16.954) (10665297,33.745) (37151361,67.000) 
};
\addplot[brown,mark=hexagon*,mark size=2.25] coordinates {
 (6561,14.027) (41553,20.022) (193185,28.004) (768609,28.002) (2935521,32.001) (10665297,34.000) (37151361,32.000) (123730257,40.000) %
};
\addlegendentry{$d$=8}
\addplot[black,domain=1e2:1e9,samples=100, thick] {ln(x)*2.5};
\addlegendentry{$\asymp \log|I|$}
\end{loglogaxis}
\end{tikzpicture}
}
\hfill
\phantom{$|$}

\hfill
\subfloat[$I:=I_{N,\mathrm{even}}^{d,0}$]{\label{fig:numerics:periodic:reco_r1l:oversampling:hc:even}
\begin{tikzpicture}
  \begin{loglogaxis}[enlargelimits=false,xmin=2.5,xmax=4e6,ymin=0.5,ymax=2e2,ytick={1,10,100},height=0.4\textwidth, width=0.94\textwidth, grid=major, xlabel={$|I_{N,\mathrm{even}}^{d,0}|$}, ylabel={$M/|I_{N,\mathrm{even}}^{d,0}|$}, font=\footnotesize,
  legend style={at={(0.5,1.05)}, anchor=south,legend columns=3,legend cell align=left, font=\footnotesize, %
  },]
\addplot[forget plot,dashed,mark options={solid},blue,mark=o,mark size=2.5] coordinates {
 (5,1.000) (13,1.000) (29,1.414) (65,2.231) (145,3.759) (329,6.422) (733,11.352) (1633,20.224) (3605,36.501) (7913,66.386) (17217,121.926)
};
\addplot[blue,mark=*,mark size=2.5] coordinates {
 (5,2.200) (13,2.231) (29,2.034) (65,2.015) (145,8.421) (329,6.052) (733,10.105) (1633,14.162) (3605,12.024) (7913,12.042) (17217,18.023)
};
\addlegendentry{$d$=2}
\addplot[forget plot,dashed,mark options={solid},red,mark=square,mark size=2] coordinates {
 (7,1.000) (25,1.160) (69,1.406) (177,2.232) (441,3.902) (1097,4.705) (2693,8.009) (6529,13.081) (15645,22.960) (37025,37.369) (86593,62.548)
};
\addplot[red,mark=square*,mark size=2] coordinates {
 (7,4.143) (25,2.120) (69,6.130) (177,10.322) (441,6.007) (1097,14.181) (2693,12.035) (6529,20.089) (15645,22.042) (37025,18.012) (86593,18.004)
};
\addlegendentry{$d$=3}
\addplot[black,domain=1e3:1e5,samples=100,dotted,very thick] {x/(ln(x))^2};
\addlegendentry{$\asymp |I|/\log^2|I|$}
\addplot[forget plot,dashed,mark options={solid},cyan,mark=triangle,mark size=2.5] coordinates {
 (9,1.000) (41,1.195) (137,1.876) (401,3.025) (1105,5.262) (2977,7.234) (7897,14.058) (20609,23.578) (52953,44.447) (133905,83.259) (333457,124.757)
};
\addplot[cyan,mark=triangle*,mark size=2.5] coordinates {
 (9,6.333) (41,6.512) (137,10.489) (401,12.257) (1105,10.090) (2977,16.137) (7897,16.039) (20609,18.012) (52953,20.009) (133905,26.006) (333457,24.003)
};
\addlegendentry{$d$=4}
\addplot[forget plot,dashed,mark options={solid},brown,mark=hexagon,mark size=2.25] coordinates {
 (17,1.000) (145,1.759) (849,3.410) (3937,5.937) (15713,11.765) (56961,21.927) (194353,56.865) (637697,132.337) %
};
\addplot[brown,mark=hexagon*,mark size=2.25] coordinates {
 (17,2.176) (145,8.421) (849,16.211) (3937,16.042) (15713,18.044) (56961,22.010) (194353,26.006) (637697,28.003) (2034289,32.001)
};
\addlegendentry{$d$=8}
\addplot[black,domain=1e2:1e8,samples=100,thick] {ln(x)*3};
\addlegendentry{$\asymp \log|I|$}
  \end{loglogaxis}
\end{tikzpicture}
}
\hfill
\phantom{$|$}
\caption{(see also \cite[Figure~2.12]{volkmerdiss}). Oversampling factors $M/|I|$ for reconstructing single rank-1 lattices (dashed lines, unfilled markers) and reconstructing multiple rank-1 lattices (solid lines, filled markers).} %
\label{fig:numerics:periodic:reco_r1l:oversampling}
\end{figure}

Since the methods based on rank-1 lattices involve an oversampling, i.e., since they require more samples than there are degrees of freedom $|I|$, we start by investigating the oversampling factors $M/|I|$, where $M$ denotes the number of samples.
\begin{Example}\label{ex:oversampling_hc}
 In Figure~\ref{fig:numerics:periodic:reco_r1l:oversampling:hc}, we visualize the obtained oversampling factors $M/|I_{N}^{d,0}|$ of reconstructing single rank-1 lattices~$\Lambda(\boldz,M)$ for symmetric hyperbolic cross index sets~$I_{N}^{d,0}$ in dimensions $d=2,3,4,8$ by dashed lines and unfilled markers. We observe that the oversampling factors are $\leq 2$ in the beginning and that they grow for increasing expansion of the hyperbolic cross~$I_{N}^{d,0}$. Note that for many realistic problem sizes runnable on current workstations, the oversampling factors still behave moderately. For instance in 8 dimensions with refinement $N=64$, we compute $|I_{64}^{d,0}|=37\,151\,361$ Fourier coefficients from $M=2\,489\,164\,387$ samples yielding an oversampling factor of 67.
 From the theoretical considerations in~\cite{Kae2012} and~\cite{kaemmererdiss}, we know that there exist reconstructing single rank-1 lattices~$\Lambda(\boldz,M)$ such that the oversampling factors are $\mathcal{O}(|I_{N}^{d,0}|/\log^d|I_{N}^{d,0}|)$, and we additionally plot these upper bounds as dotted graphs. The observed oversampling factors approximately behave like these upper bounds. \newline
 Moreover, we display the oversampling factors for the reconstructing multiple rank-1 lattices~$\Lambda$ generated by Algorithm~\ref{alg:construct_mr1l_I_distinct_primes} as solid lines and filled markers. We observe that the initial oversampling factors are distinctly larger for small to medium cardinalities $|I_{N}^{d,0}|$ compared to the single reconstructing rank-1 lattices~$\Lambda(\boldz,M)$. However, the oversampling factors for the reconstructing multiple rank-1 lattice~$\Lambda$ grow slower and for larger cardinalities~$|I_{N}^{d,0}|$, the multiple rank-1 lattices yield smaller oversampling factors. As discussed in \cite{Kae17}, the multiple rank-1 lattices generated by Algorithm~\ref{alg:construct_mr1l_I_distinct_primes} have an oversampling factor of only $\mathcal{O}(\log |I|)$ with high probability. Additionally, we plot this upper bound in Figure~\ref{fig:numerics:periodic:reco_r1l:oversampling:hc} and we observe that the obtained oversampling factors for the multiple rank-1 lattices~$\Lambda$ seem to behave accordingly.\newline
 We remark that reconstructing multiple rank-1 lattices~$\Lambda$ have another important advantage over single rank-1 lattices. There are distinctly faster construction algorithms available, such as Algorithm~\ref{alg:construct_mr1l_I_distinct_primes}, requiring only $\mathcal{O}(|I|(d+\log|I|)\log|I|)$ arithmetic operations with high probability in contrast to $\mathcal{O}(|I|^3+d|I|^2\log|I|)$ arithmetic operations for reconstructing single rank-1 lattices.
 For instance, the generation of the reconstructing single rank-1 lattice~$\Lambda(\boldz,M)$ for $I=I_{64}^{8,0}$,
 consisting of $M=2\,489\,164\,387$ sampling nodes for $|I|=37\,151\,361$, required several days on a computer with Intel Xeon E7-4880 v2 CPU (2.50~GHz) using 12 threads, whereas Algorithm~\ref{alg:construct_mr1l_I_distinct_primes} returned a reconstructing multiple rank-1 lattice~$\Lambda$, consisting of $L=16$ single rank-1 lattices with $|\Lambda|=1\,188\,846\,719$ sampling nodes, within less than 1~hour using 1~thread.
 We remark that the numbers $L$ of single rank-1 lattices
 that are combined to achieve reconstructing multiple rank-1 lattices
 are in the range $\{1,\ldots,20\}$ in all numerical experiments.
\end{Example}

Next, we compare the resulting sampling errors when using the reconstructing single and reconstructing multiple rank-1 lattices as sampling sets as discussed above.

\begin{Example} \label{ex:G34_rel_sampl_err}
 We sample the test functions $G_{3,4}^d$ in dimensions $d\in\{2,\ldots,8\}$ along reconstructing single and reconstructing multiple rank-1 lattices for symmetric hyperbolic cross index sets $I_N^{d,0}$ of various refinements $N\in\N$ using the sampling operators~\eqref{eqn:sampling_operator} and~\eqref{eqn:sampling_operator_mr1l}.
 We determine the sampling errors in the relative $\mathcal{A}(\T^d)$ norm and relative $L_2(\T^d)$ norm. The corresponding results are depicted in Figures~\ref{fig:numerics:periodic:G34:rel_sampl_err_A:M} and~\ref{fig:numerics:periodic:G34:rel_sampl_err_L2:M}, respectively, where the results for single rank-1 lattices are plotted as dashed lines with unfilled markers and the results for multiple rank-1 lattices as solid lines with filled markers.
 We observe that the sampling errors for single rank-1 lattices decrease slower in general than for multiple rank-1 lattices.
 When using the same hyperbolic cross frequency index sets $I=I_N^{d,0}$, the single rank-1 lattices yield slightly smaller error values than the multiple ones.
 Correspondingly, since the oversampling factors of the single rank-1 lattices are lower for smaller index sets~$I$, cf.\ Example~\ref{ex:oversampling_hc}, the single rank-1 lattices perform better in these cases when considering the sampling error with respect to the numbers of samples. Once the frequency index sets~$I$ become larger, the multiple rank-1 lattices perform better.
 Again, we also stress on the fact that the construction of reconstructing single rank-1 lattices may require much more time compared to multiple rank-1 lattices, see also the discussion at the end of Example~\ref{ex:oversampling_hc}.
 Since the test functions $G_{3,4}^d\in\mathcal{A}^{0,3-\epsilon}(\T^d)$ and $G_{3,4}^d\in\mathcal{H}^{0,\frac{7}{2}-\epsilon}(\T^d)$, $\epsilon > 0$,
 we expect the relative sampling errors for single rank-1 lattices in Figure~\ref{fig:numerics:periodic:G34:rel_sampl_err_A:M} and~\ref{fig:numerics:periodic:G34:rel_sampl_err_L2:M} to almost decay like $\sim M^{-\frac{3-\epsilon}{2}} \, (\log M)^{\frac{d-2}{2}(3-\epsilon)}$ and $\sim M^{-\frac{3.5-\epsilon}{2}} \, (\log M)^{\frac{d-2}{2}(3.5-\epsilon)+\frac{d-1}{2}}$, cf.\ \cite[Corollaries~2.40 and~2.44]{volkmerdiss}, respectively. These upper bounds are plotted in Figure~\ref{fig:numerics:periodic:G34:rel_sampl_err} as dotted graphs for dimensions $d=2,3,4$ and we observe that the relative $\mathcal{A}(\T^d)$ and $L_2(\T^d)$ sampling errors for single rank-1 lattices approximately behave like these bounds.
 For multiple rank-1 lattices, we expect from Corollaries~\ref{cor:sampl_error_mr1l:Linf:Aabg} and~\ref{cor:approx_error:Hrtg:Habg} that the relative $\mathcal{A}(\T^d)$ and $L_2(\T^d)$ sampling errors should decay like $\sim M^{-(3-\epsilon)} \, (\log M)^{(3-\epsilon)d+1}$ and $\sim M^{-(3.5-\epsilon-\lambda)} \, (\log M)^{(3.5-\epsilon-\lambda)d+1}$, respectively, where $\lambda>1/2$. We visualize these upper bounds by solid lines without markers in Figure~\ref{fig:numerics:periodic:G34:rel_sampl_err} for $d=2,3,4$ and we observe that the obtained errors approximately behave like these bounds suggest or slightly better. The error plots for dimensions greater than three still suffer from pre-asymptotic behavior. For that reason, we have left out the corresponding plots of the asymptotic bounds. Moreover, for better recognizability, we omitted the plots for $d=7$ which show a behavior in-between those of $d=6$ and $d=8$.
\end{Example}

\begin{figure}[htb]
\centering
\subfloat[$\Vert G_{3,4}^d - S^\Lambda_{I} G_{3,4}^d\vert \mathcal{A}(\T^d)\Vert / \Vert G_{3,4}^d\vert \mathcal{A}(\T^d)\Vert$]{\label{fig:numerics:periodic:G34:rel_sampl_err_A:M}
\begin{tikzpicture}[baseline]
  \begin{loglogaxis}[enlargelimits=false,xmin=5e0,xmax=1e10,ymin=1e-9,ymax=3e0,ytick={1e-8,1e-6,1e-4,1e-2,1},height=0.4\textwidth, width=0.95\textwidth, grid=major, xlabel={number of samples $M$}, ylabel={rel\_{}sampl\_{}err\_A}, font=\footnotesize,
  legend style={at={(0.5,1.05)}, anchor=south,legend columns=4,legend cell align=left, font=\footnotesize, 
  },
  xminorticks=false,yminorticks=false
  ]
\addplot[forget plot,dashed,mark options={solid},blue,mark=o,mark size=2.5] coordinates {
 (9,4.832e-01) (23,2.496e-01) (58,5.852e-02) (163,1.009e-02) (579,9.979e-04) (2179,1.093e-04) (8451,1.325e-05) (33283,1.650e-06) (132099,2.106e-07) (526339,2.693e-08) (2101251,3.478e-09)
};
\addplot[blue,mark=*,mark size=2.5] coordinates {
 (17,4.237e-01) (41,1.809e-01) (197,4.466e-02) (687,1.209e-02) (1087,1.012e-03) (3651,1.139e-04) (13873,1.584e-05) (37243,1.937e-06) (69001,2.499e-07) (212659,3.341e-08) (397975,4.244e-09) %
};
\addlegendentry{$d$=2}
\addplot[forget plot,dashed,mark options={solid},red,mark=square,mark size=2] coordinates {
 (27,6.533e-01) (105,3.481e-01) (343,1.103e-01) (1035,2.918e-02) (3628,5.124e-03) (11525,8.075e-04) (47463,1.042e-04) (176603,1.313e-05) (753249,1.716e-06) (2773801,2.268e-07) (11121676,3.029e-08) (47533463,4.044e-09)
};
\addplot[red,mark=square*,mark size=2] coordinates {
 (53,4.773e-01) (859,2.903e-01) (1365,9.820e-02) (3579,3.474e-02) (28679,7.009e-03) (56525,1.114e-03) (162143,1.369e-04) (497753,1.650e-05) (1205421,2.225e-06) (3437839,3.001e-07) (7393305,3.969e-08) (17150815,5.286e-09)
};
\addlegendentry{$d$=3}
\addplot[forget plot,dashed,mark options={solid},cyan,mark=triangle,mark size=2.5] coordinates {
 (81,7.766e-01) (479,4.318e-01) (1911,1.601e-01) (5727,5.763e-02) (21944,1.401e-02) (106703,2.866e-03) (475829,5.020e-04) (2244100,8.109e-05) (10561497,1.141e-05) (39632648,1.578e-06) (176847961,2.184e-07)
};
\addplot[cyan,mark=triangle*,mark size=2.5] coordinates {
 (679,6.904e-01) (3625,4.103e-01) (7607,1.618e-01) (22259,6.782e-02) (130115,1.727e-02) (453711,3.817e-03) (1238275,6.342e-04) (4268363,1.075e-04) (10229923,1.530e-05) (30306295,2.103e-06) (86364359,2.916e-07) (212169309,4.120e-08)
};
\addlegendentry{$d$=4}
\addplot[black,domain=1e5:1e9,samples=100,dotted,very thick] {30*x^(-3/2)*ln(x)^((2-2)*3/2)};
\addlegendentry{theor. single}
\addplot[forget plot,black,domain=1e6:1e8,samples=100,dotted,very thick] {50*x^(-3/2)*ln(x)^((3-2)*3/2)};
\addplot[forget plot,black,domain=1e7:5e8,samples=100,dotted,very thick] {300*x^(-3/2)*ln(x)^((4-2)*3/2)};
\addplot[forget plot,dashed,mark options={solid},darkgreen,mark=diamond,mark size=2.5] coordinates {
 (243,8.653e-01) (2185,5.053e-01) (10579,2.152e-01) (33769,8.979e-02) (169230,2.750e-02) (785309,7.106e-03) (3752318,1.571e-03) (20645268,3.150e-04) (136178715,5.455e-05)
};
\addplot[darkgreen,mark=diamond*,mark size=2.5] coordinates {
 (1475,8.385e-01) (8481,4.903e-01) (45149,2.376e-01) (168623,1.080e-01) (840967,3.548e-02) (2308173,9.560e-03) (7443841,2.161e-03) (28751753,4.353e-04) (74060555,7.532e-05) (199219495,1.197e-05)
};
\addlegendentry{$d$=5}
\addplot[forget plot,dashed,mark options={solid,rotate=180},darkorange,mark=triangle,mark size=2.5] coordinates {
 (729,9.295e-01) (9967,5.710e-01) (57897,2.720e-01) (191808,1.295e-01) (1105193,4.626e-02) (6897012,1.417e-02) (31829977,3.755e-03) (192757285,8.938e-04) (1400567254,1.863e-04)
};
\addplot[darkorange,mark options={rotate=180},darkorange,mark=triangle*,mark size=2.5] coordinates {
 (5891,9.482e-01) (65959,6.235e-01) (200995,3.187e-01) (995579,1.611e-01) (5077743,6.002e-02) (14234871,1.875e-02) (47878105,5.102e-03) (133583559,1.206e-03) (569137025,2.578e-04) (1527406155,4.898e-05)
};
\addlegendentry{$d$=6}
\addplot[forget plot,dashed,mark options={solid},brown,mark=hexagon,mark size=2.25] coordinates {
 (6561,1.007e+00) (207391,6.804e-01) (1259193,3.964e-01) (6027975,2.170e-01) (49768670,9.905e-02) (359896131,3.887e-02) (2489164387,1.330e-2)
};
\addplot[brown,mark=hexagon*,mark size=2.25] coordinates {
 (92033,9.830e-01) (831963,7.503e-01) (5409995,4.795e-01) (21522733,2.742e-01) (93939469,1.322e-01) (362622757,5.217e-02) (1188846719,1.795e-02) (4949213695,5.617e-03)
};
\addlegendentry{$d$=8}
\addplot[black,domain=1e4:1e6,samples=100,thick] {1/1.5*x^(-3)*ln(x)^(2*3+1)};
\addlegendentry{theor. mult.}
\addplot[forget plot,black,domain=8e5:9e6,samples=100,thick] {1.8*x^(-3)*ln(x)^(3*3+1)};
\addplot[forget plot,black,domain=1e7:1.5e8,samples=100,thick] {1.2*x^(-3)*ln(x)^(4*3+1)};
\legend{};
  \end{loglogaxis}
\end{tikzpicture}
}
\vspace{1em}
\subfloat[$\Vert G_{3,4}^d - S^\Lambda_{I} G_{3,4}^d\vert L_2(\T^d)\Vert / \Vert G_{3,4}^d\vert L_2(\T^d)\Vert$]{\label{fig:numerics:periodic:G34:rel_sampl_err_L2:M}
\begin{tikzpicture}[baseline]
  \begin{loglogaxis}[enlargelimits=false,xmin=5e0,xmax=1e10,ymin=1e-11,ymax=3e0,ytick={1e-10,1e-8,1e-6,1e-4,1e-2,1},height=0.4\textwidth, width=0.95\textwidth, grid=major, xlabel={number of samples $M$}, ylabel={rel\_{}sampl\_{}err\_L2}, font=\footnotesize,
  legend style={legend columns=4,legend cell align=left, font=\footnotesize},
  legend to name=namedcgsteps,    
  xminorticks=false,yminorticks=false
  ]
\addplot[forget plot,dashed,mark options={solid},blue,mark=o,mark size=2.5] coordinates {
(9,2.165e-01) (23,1.044e-01) (58,1.679e-02) (163,3.606e-03) (579,2.028e-04) (2179,1.295e-05) (8451,1.129e-06) (33283,1.005e-07) (132099,8.969e-09) (526339,7.954e-10) (2101251,7.163e-11)
};
\addplot[blue,mark=*,mark size=2.5] coordinates {
 (17,1.828e-01) (41,8.431e-02) (197,1.346e-02) (687,3.796e-03) (1087,2.027e-04) (3651,1.292e-05) (13873,1.174e-06) (37243,1.038e-07) (69001,9.232e-09) (212659,8.262e-10) (397975,7.405e-11) %
};
\addlegendentry{$d$=2}
\addplot[forget plot,dashed,mark options={solid},red,mark=square,mark size=2] coordinates {
 (27,2.675e-01) (105,1.227e-01) (343,2.320e-02) (1035,6.338e-03) (3628,7.317e-04) (11525,9.875e-05) (47463,7.159e-06) (176603,4.989e-07) (753249,4.829e-08) (2773801,3.943e-09) (11121676,3.608e-10) (47533463,3.262e-11)
};
\addplot[red,mark=square*,mark size=2] coordinates {
 (53,1.995e-01) (859,1.069e-01) (1365,2.033e-02) (3579,6.723e-03) (28679,7.793e-04) (56525,1.069e-04) (162143,7.488e-06) (497753,5.136e-07) (1205421,5.053e-08) (3437839,4.078e-09) (7393305,3.727e-10) (17150815,3.359e-11)
};
\addlegendentry{$d$=3}
\addplot[forget plot,dashed,mark options={solid},cyan,mark=triangle,mark size=2.5] coordinates {
 (81,3.086e-01) (479,1.380e-01) (1911,2.865e-02) (5727,9.359e-03) (21944,1.387e-03) (106703,2.153e-04) (475829,2.605e-05) (2244100,3.493e-06) (10561497,2.529e-07) (39632648,1.937e-08) (176847961,1.699e-09)
};
\addplot[cyan,mark=triangle*,mark size=2.5] coordinates {
 (679,2.486e-01) (3625,1.271e-01) (7607,2.671e-02) (22259,9.789e-03) (130115,1.417e-03) (453711,2.237e-04) (1238275,2.669e-05) (4268363,3.589e-06) (10229923,2.602e-07) (30306295,1.990e-08) (86364359,1.742e-09)
};
\addlegendentry{$d$=4}
\addplot[black,domain=1e5:5e6,samples=100,dotted,very thick] {8*x^(-3.5/2)*ln(x)^((2-2)*3.5/2+(2-1)/2)};
\addlegendentry{theor. single} %
\addplot[forget plot,black,domain=3e6:1e9,samples=100,dotted,very thick] {2*x^(-3.5/2)*ln(x)^((3-2)*3.5/2+(3-1)/2)};
\addplot[forget plot,black,domain=1e7:1e9,samples=100,dotted,very thick] {1.5*x^(-3.5/2)*ln(x)^((4-2)*3.5/2+(4-1)/2)};
\addplot[forget plot,dashed,mark options={solid},darkgreen,mark=diamond,mark size=2.5] coordinates {
 (243,3.440e-01) (2185,1.519e-01) (10579,3.426e-02) (33769,1.179e-02) (169230,2.131e-03) (785309,3.693e-04) (3752318,5.499e-05) (20645268,8.578e-06) (136178715,9.776e-07)
};
\addplot[darkgreen,mark=diamond*,mark size=2.5] coordinates {
 (1475,3.077e-01) (8481,1.404e-01) (45149,3.326e-02) (168623,1.231e-02) (840967,2.194e-03) (2308173,3.829e-04) (7443841,5.707e-05) (28751753,8.846e-06) (74060555,1.007e-06) (199219495,1.169e-07)
};
\addlegendentry{$d$=5}
\addplot[forget plot,dashed,mark options={solid,rotate=180},darkorange,mark=triangle,mark size=2.5] coordinates {
 (729,3.752e-01) (9967,1.649e-01) (57897,4.000e-02) (191808,1.458e-02) (1105193,2.957e-03) (6897012,5.579e-04) (31829977,9.469e-05) (192757285,1.627e-05) (1400567254,2.272e-06)
};
\addplot[darkorange,mark options={rotate=180},darkorange,mark=triangle*,mark size=2.5] coordinates {
 (5891,3.206e-01) (65959,1.591e-01) (200995,3.966e-02) (995579,1.511e-02) (5077743,3.033e-03) (14234871,5.740e-04) (47878105,9.751e-05) (133583559,1.672e-05) (569137025,2.324e-06) (1527406155,3.136e-07)
};
\addlegendentry{$d$=6}
\addplot[forget plot,dashed,mark options={solid},brown,mark=hexagon,mark size=2.25] coordinates {
 (6561,4.289e-01) (207391,1.886e-01) (1259193,5.225e-02) (6027975,2.023e-02) (49768670,4.879e-03) (359896131,1.046e-03)  (2489164387,2.086e-04)
};
\addplot[brown,mark=hexagon*,mark size=2.25] coordinates {
 (92033,3.407e-01) (831963,1.800e-01) (5409995,5.182e-02) (21522733,2.076e-02) (93939469,5.015e-03) (362622757,1.071e-03) (1188846719,2.137e-04) (4949213695,4.184e-05)
};
\addlegendentry{$d$=8}
\addplot[black,domain=1e4:1e6,samples=100,thick] {1/100*x^(-3)*ln(x)^(2*3+1)};
\addlegendentry{theor. mult.}
\addplot[forget plot,black,domain=3e5:1e8,samples=100,thick] {1/50*x^(-3)*ln(x)^(3*3+1)};
\addplot[forget plot,black,domain=3e6:1e8,samples=100,thick] {1/100*x^(-3)*ln(x)^(4*3+1)};
  \end{loglogaxis}
\end{tikzpicture}
}
\\[1em]
\pgfplotslegendfromname{namedcgsteps}
\caption{(see also \cite[Figure~2.10]{volkmerdiss}). Relative $\mathcal{A}(\T^d)$ and
 $L_2(\T^d)$ sampling errors 
  for $G_{3,4}^d$ with respect to the number of sampling nodes $M$ for reconstructing single rank-1 lattices (dashed lines, unfilled markers) and reconstructing multiple rank-1 lattices (solid lines, filled markers),
  when using the frequency index sets $I:=I_{N}^{d,0}$.
  }\label{fig:numerics:periodic:G34:rel_sampl_err}
\end{figure}

Additionally,
we consider the tensor-product test functions
$G_3^d\colon\T^d\rightarrow\C$ from~\cite{KaPoVo13},
$G_3^d(\boldx):=\prod_{s=1}^{d} g_3(x_s)$,
where the one-dimensional function $g_3\colon\T\rightarrow\C$ is defined by
$$
g_3(x):= 4 \, \sqrt{\frac{3\pi}{207 \pi - 256}} \left(2+\sgn((x\bmod 1)-1/2) \, \sin(2\pi x)^3\right)
$$
and
$\Vert G_3^d \vert L^2(\T^d) \Vert = 1$.
We have
$G_3^d\in\mathcal{A}^{0,3-\epsilon}(\T^d)$ and $G_3^d\in\mathcal{H}^{0,\frac{7}{2}-\epsilon}(\T^d)$, $\epsilon > 0$,
as well as $\Vert G_3^d \vert \mathcal{A}(\T^d) \Vert = \left(\frac{8 (4 + 15 \pi)}{5 \sqrt{3 \pi} (207 \pi - 256)}\right)^d \approx (1.34181)^d$.
Since the Fourier coefficients $(\widehat{g_3})_k$ of $g_3$ are zero for odd frequencies $k\in(2\Z+1)$,
we consider hyperbolic cross index sets ``with holes'' $I=I_{N,\mathrm{even}}^{d,0}:=I_N^{d,0} \cap (2\Z)^d$,
which consist of distinctly less frequencies compared to the index sets $I_N^{d,0}$.

\begin{Example}\label{ex:oversampling_hc_even}
 In Figure~\ref{fig:numerics:periodic:reco_r1l:oversampling:hc:even}, we depict the oversampling factors $M/|I_{N,\mathrm{even}}^{d,0}|$ of reconstructing single rank-1 lattices~$\Lambda(\boldz,M)$ for symmetric hyperbolic cross index sets ``with holes'' $I=I_{N,\mathrm{even}}^{d,0}$ generated by the implementation \cite[\texttt{genlattice\_cbc\_incr\_bisect}]{Vo_taylorR1Lnfft} of \cite[Algorithm~3.7]{kaemmererdiss}, see also \cite[Table~2.3]{volkmerdiss}, by dashed lines and unfilled markers.
 Additionally, we visualize the oversampling factors for the reconstructing multiple rank-1 lattices~$\Lambda$ generated by Algorithm~\ref{alg:construct_mr1l_I_distinct_primes} as solid lines and filled markers.
 We observe almost the same behavior as in Example~\ref{ex:oversampling_hc} for symmetric hyperbolic cross index sets $I_{N}^{d,0}$.
\end{Example}

\begin{Example} \label{ex:G3_rel_sampl_err}
 Now, we sample the test functions $G_3^d$ along reconstructing single and multiple rank-1 lattices for symmetric hyperbolic cross index sets ``with holes'' $I_{N,\mathrm{even}}^{d,0}$ of various refinements $N\in\N$ using the sampling operators~\eqref{eqn:sampling_operator} and~\eqref{eqn:sampling_operator_mr1l}.
 As in Example~\ref{ex:G34_rel_sampl_err}, we determine the sampling errors in the relative $\mathcal{A}(\T^d)$ norm and relative $L_2(\T^d)$ norm. The corresponding results are depicted in Figures~\ref{fig:numerics:periodic:G3:rel_sampl_err_A:M} and~\ref{fig:numerics:periodic:G3:rel_sampl_err_L2:M}, respectively, for dimensions $d=2,3,4,6,8$, where the results for single rank-1 lattices are plotted as dashed lines with unfilled markers and the results for multiple rank-1 lattices as solid lines with filled markers.
 In principle, we observe almost the same behavior as in Example~\ref{ex:G34_rel_sampl_err}, since the test functions are also in $\mathcal{A}^{0,3-\epsilon}(\T^d)$ and $\mathcal{H}^{0,\frac{7}{2}-\epsilon}(\T^d)$, $\epsilon > 0$.
 Again, we plot the graphs of the theoretical upper bounds as in Example~\ref{ex:G34_rel_sampl_err}.
 When using multiple rank-1 lattices and large refinements~$N$, we notice that the relative $L_2(\T^d)$ sampling errors seem to decay slightly faster than the theoretical upper bounds from Corollary~\ref{cor:approx_error:Hrtg:Habg} for dimensions $d=2,3,4$. We suspect that the additive term $\lambda>1/2$ may not occur.
 \newline
 Additionally, we visualize the relative sampling errors in higher dimensions up to $d=20$ in Figure~\ref{fig:numerics:periodic:G3:rel_sampl_err:M:hd}, which are computable in practice due to the sparser structure of the hyperbolic cross index sets ``with holes'' $I=I_{N,\mathrm{even}}^{d,0}$. Since the oversampling factors of reconstructing single rank-1 lattices are lower for smaller frequency index sets~$I$, cf.\ Example~\ref{ex:oversampling_hc_even} and Figure~\ref{fig:numerics:periodic:reco_r1l:oversampling:hc:even}, using single rank-1 lattices requires less samples for comparable sampling errors, whereas multiple rank-1 lattices are better suited for larger frequency index sets~$I$. However, we emphasize once more that building the reconstructing multiple rank-1 lattices typically requires distinctly less runtime, see also Example~\ref{ex:oversampling_hc}.
\end{Example}

\begin{figure}[htb]
\centering
\subfloat[$\Vert G_3^d - S^\Lambda_{I} G_3^d\vert \mathcal{A}(\T^d)\Vert / \Vert G_3^d\vert \mathcal{A}(\T^d)\Vert$]{\label{fig:numerics:periodic:G3:rel_sampl_err_A:M}
\begin{tikzpicture}
  \begin{loglogaxis}[enlargelimits=false,xmin=4,xmax=1e9,ymin=4e-9,ymax=3e0,ytick={1e-8,1e-6,1e-4,1e-2,1},height=0.4\textwidth, width=0.95\textwidth, grid=major, xlabel={number of samples $M$}, ylabel={rel\_{}sampl\_{}err\_A}, font=\footnotesize,
  legend style={at={(0.375,1.05)}, anchor=south,legend columns=4,legend cell align=left, font=\footnotesize, %
  },
  ]
\addplot[forget plot,dashed,mark options={solid},blue,mark=o,mark size=2.5] coordinates {
 (5,1.504e-01) (13,5.298e-02) (41,6.735e-03) (145,1.124e-03) (545,1.837e-04) (2113,2.639e-05) (8321,3.830e-06) (33025,5.302e-07) (131585,7.248e-08) (525313,9.751e-09)
};
\addplot[blue,mark=*,mark size=2.5] coordinates {
 (11,1.404e-01) (29,3.980e-02) (59,6.045e-03) (131,1.193e-03) (1221,1.993e-04) (1991,2.936e-05) (7407,4.419e-06) (23127,6.512e-07) (43347,8.431e-08) (95289,1.198e-08) %
};
\addlegendentry{$d$=2}
\addplot[forget plot,dashed,mark options={solid},red,mark=square,mark size=2] coordinates {
 (7,2.908e-01) (29,1.166e-01) (97,2.476e-02) (395,4.382e-03) (1721,8.530e-04) (5161,1.585e-04) (21569,2.651e-05) (85405,4.160e-06) (359213,6.266e-07) (1383595,9.214e-08) (5416219,1.330e-08)
};
\addplot[red,mark=square*,mark size=2] coordinates {
 (29,3.278e-01) (53,8.239e-02) (423,2.397e-02) (1827,5.013e-03) (2649,1.013e-03) (15557,2.002e-04) (32411,3.320e-05) (131163,5.013e-06) (344841,7.596e-07) (666907,1.156e-07) (1559059,1.668e-08)
};
\addlegendentry{$d$=3}
\addplot[forget plot,dashed,mark options={solid},cyan,mark=triangle,mark size=2.5] coordinates {
 (9,4.551e-01) (49,1.681e-01) (257,5.163e-02) (1213,1.209e-02) (5815,2.855e-03) (21535,6.215e-04) (111015,1.242e-04) (485913,2.275e-05) (2353599,3.866e-06) (11148851,6.256e-07) (41601005,9.821e-08)
};
\addplot[cyan,mark=triangle*,mark size=2.5] coordinates {
 (57,4.249e-01) (267,1.647e-01) (1437,5.102e-02) (4915,1.375e-02) (11149,3.446e-03) (48041,7.574e-04) (126659,1.507e-04) (371217,2.834e-05) (1059531,4.843e-06) (3482387,7.772e-07) (8004093,1.209e-07)
};
\addlegendentry{$d$=4}
\addplot[black,domain=1e4:8e5,samples=100,dotted,very thick] {10*x^(-3/2)*ln(x)^((2-2)*3/2)};
\addlegendentry{theor. single} %
\addplot[forget plot,black,domain=4e5:1e7,samples=100,dotted,very thick] {8*x^(-3/2)*ln(x)^((3-2)*3/2)};
\addplot[forget plot,black,domain=1e6:1e8,samples=100,dotted,very thick] {15*x^(-3/2)*ln(x)^((4-2)*3/2)};
\addplot[forget plot,dashed,mark options={solid,rotate=180},darkorange,mark=triangle,mark size=2.5] coordinates { %
 (13,7.515e-01) (137,3.472e-01) (983,1.492e-01) (6905,4.944e-02) (34117,1.597e-02) (226951,4.469e-03) (1373325,1.178e-03) (8145033,2.859e-04) (50770301,6.387e-05) (293168219,1.314e-05)
};
\addplot[darkorange,mark options={rotate=180},darkorange,mark=triangle*,mark size=2.5] coordinates {
 (95,7.121e-01) (721,3.524e-01) (5669,1.507e-01) (14689,5.486e-02) (78061,1.751e-02) (274677,5.257e-03) (1198585,1.392e-03) (3262425,3.389e-04) (10211917,7.676e-05) (31145671,1.630e-05)
};
\addlegendentry{$d$=6}
\addplot[forget plot,dashed,mark options={solid},brown,mark=hexagon,mark size=2.25] coordinates {
 (17,9.667e-01) (255,5.241e-01) (2895,2.770e-01) (23375,1.231e-01) (184859,4.868e-02) (1248979,1.742e-02) (11051805,5.657e-03) (84391053,1.706e-03) (600266399,4.745e-04)
};
\addplot[brown,mark=hexagon*,mark size=2.25] coordinates {
 (37,8.009e-01) (1221,5.098e-01) (13763,2.762e-01) (63159,1.324e-01) (283521,5.289e-02) (1253697,1.934e-02) (5054425,6.474e-03) (17857181,1.967e-03) (65099159,5.580e-04) (202803915,1.461e-04)
};
\addlegendentry{$d$=8}
\addlegendimage{empty legend}
\addlegendentry{}
\addplot[black,domain=4e3:1e6,samples=100,thick] {0.1*x^(-3)*ln(x)^(2*3+1)};
\addlegendentry{theor. mult.}
\addplot[forget plot,black,domain=1e5:1e6,samples=100,thick] {0.05*x^(-3)*ln(x)^(3*3+1)};
\addplot[forget plot,black,domain=4e5:7e6,samples=100,thick] {0.003*x^(-3)*ln(x)^(4*3+1)};
\legend{};
  \end{loglogaxis}
\end{tikzpicture}
}
\vspace{1em}
\subfloat[$\Vert G_3^d - S^\Lambda_{I} G_3^d\vert L_2(\T^d)\Vert / \Vert G_3^d\vert L_2(\T^d)\Vert$]{\label{fig:numerics:periodic:G3:rel_sampl_err_L2:M}
\begin{tikzpicture}
  \begin{loglogaxis}[enlargelimits=false,xmin=3,xmax=1e9,ymin=7e-11,ymax=2e0,ytick={1e-8,1e-6,1e-4,1e-2,1},height=0.4\textwidth, width=0.95\textwidth, grid=major, xlabel={number of samples $M$}, ylabel={rel\_{}sampl\_{}err\_L2}, font=\footnotesize,
  legend style={at={(0.375,1.05)}, anchor=south,legend columns=4,legend cell align=left, font=\footnotesize, %
  },
  legend to name=samperrG3dM,
  ]
\addplot[forget plot,dashed,mark options={solid},blue,mark=o,mark size=2.5] coordinates {
 (5,6.973e-02) (13,1.939e-02) (41,1.854e-03) (145,2.324e-04) (545,2.846e-05) (2113,2.789e-06) (8321,2.796e-07) (33025,2.625e-08) (131585,2.437e-09) (525313,2.226e-10)
};
\addplot[blue,mark=*,mark size=2.5] coordinates {
 (11,6.759e-02) (29,1.386e-02) (59,1.688e-03) (131,2.390e-04) (1221,2.929e-05) (1991,2.856e-06) (7407,2.912e-07) (23127,2.758e-08) (43347,2.519e-09) (95289,2.374e-10) %
};
\addlegendentry{$d$=2}
\addplot[forget plot,dashed,mark options={solid},red,mark=square,mark size=2] coordinates {
 (7,1.091e-01) (29,3.408e-02) (97,4.736e-03) (395,6.420e-04) (1721,9.394e-05) (5161,1.271e-05) (21569,1.514e-06) (85405,1.527e-07) (359213,1.525e-08) (1383595,1.464e-09) (5416219,1.385e-10)
};
\addplot[red,mark=square*,mark size=2] coordinates {
 (29,1.387e-01) (53,2.273e-02) (423,4.627e-03) (1827,6.747e-04) (2649,9.872e-05) (15557,1.367e-05) (32411,1.604e-06) (131163,1.577e-07) (344841,1.574e-08) (666907,1.533e-09) (1559059,1.448e-10)
};
\addlegendentry{$d$=3}
\addplot[forget plot,dashed,mark options={solid},cyan,mark=triangle,mark size=2.5] coordinates {
 (9,1.685e-01) (49,4.193e-02) (257,8.063e-03) (1213,1.353e-03) (5815,2.294e-04) (21535,3.593e-05) (111015,5.188e-06) (485913,6.736e-07) (2353599,8.093e-08) (11148851,7.810e-09) (41601005,7.895e-10)
};
\addplot[cyan,mark=triangle*,mark size=2.5] coordinates {
  (57,1.514e-01) (267,3.782e-02) (1437,7.668e-03) (4915,1.418e-03) (11149,2.459e-04) (48041,3.775e-05) (126659,5.398e-06) (371217,7.030e-07) (1059531,8.363e-08) (3482387,8.088e-09) (8004093,8.153e-10)
};
\addlegendentry{$d$=4}
\addplot[black,domain=1e4:9e5,samples=100,dotted,very thick] {1.5*x^(-3.5/2)*ln(x)^((2-2)*3.5/2+(2-1)/2)};
\addlegendentry{theor. single} %
\addplot[forget plot,black,domain=1e5:8e6,samples=100,dotted,very thick] {0.15*x^(-3.5/2)*ln(x)^((3-2)*3.5/2+(3-1)/2)};
\addplot[forget plot,black,domain=1e6:1e8,samples=100,dotted,very thick] {0.05*x^(-3.5/2)*ln(x)^((4-2)*3.5/2+(4-1)/2)};
\addplot[forget plot,dashed,mark options={solid,rotate=180},darkorange,mark=triangle,mark size=2.5] coordinates {
 (13,3.643e-01) (137,7.194e-02) (983,1.800e-02) (6905,3.818e-03) (34117,8.144e-04) (226951,1.552e-04) (1373325,2.813e-05) (8145033,4.750e-06) (50770301,7.462e-07) (293168219,1.115e-07)
};
\addplot[darkorange,mark options={rotate=180},darkorange,mark=triangle*,mark size=2.5] coordinates {
 (95,3.179e-01) (721,7.270e-02) (5669,1.790e-02) (14689,4.031e-03) (78061,8.317e-04) (274677,1.632e-04) (1198585,2.913e-05) (3262425,4.891e-06) (10211917,7.668e-07) (31145671,1.145e-07)
};
\addlegendentry{$d$=6}
\addplot[forget plot,dashed,mark options={solid},brown,mark=hexagon,mark size=2.25] coordinates {
 (17,6.810e-01) (255,1.087e-01) (2895,3.003e-02) (23375,7.830e-03) (184859,1.909e-03) (1248979,4.331e-04) (11051805,9.125e-05) (84391053,1.814e-05) (600266399,3.389e-06)
};
\addplot[brown,mark=hexagon*,mark size=2.25] coordinates {
 (37,3.336e-01) (1221,1.016e-01) (13763,2.974e-02) (63159,8.251e-03) (283521,1.976e-03) (1253697,4.470e-04) (5054425,9.463e-05) (17857181,1.868e-05) (65099159,3.485e-06) (202803915,6.122e-07)
};
\addlegendentry{$d$=8}
\addlegendimage{empty legend}
\addlegendentry{}
\addplot[black,domain=1e4:1e6,samples=100,thick] {0.002*x^(-3)*ln(x)^(2*3+1)};
\addlegendentry{theor. mult.}
\addplot[forget plot,black,domain=1e5:1e6,samples=100,thick] {0.0009*x^(-3)*ln(x)^(3*3+1)};
\addplot[forget plot,black,domain=4e5:9e6,samples=100,thick] {0.00005*x^(-3)*ln(x)^(4*3+1)};
  \end{loglogaxis}
\end{tikzpicture}
}
\\[1em]
\pgfplotslegendfromname{samperrG3dM}
\caption{(see also \cite[Figure~2.14]{volkmerdiss}). Relative $\mathcal{A}(\T^d)$
and $L_2(\T^d)$ sampling errors
for $G_3^d$ with respect to the number of sampling nodes~$M$ for reconstructing single rank-1 lattices (dashed lines, unfilled markers) and reconstructing multiple rank-1 lattices (solid lines, filled markers),
when using the frequency index sets $I:=I_{N,\mathrm{even}}^{d,0}$.
}\label{fig:numerics:periodic:G3:rel_sampl_err:M}
\end{figure}

\begin{figure}[htb]
\centering
\subfloat[$\Vert G_3^d - S^\Lambda_{I} G_3^d\vert \mathcal{A}(\T^d)\Vert / \Vert G_3^d\vert \mathcal{A}(\T^d)\Vert$]{\label{fig:numerics:periodic:G3:rel_sampl_err_A:M:hd}
\begin{tikzpicture}
  \begin{loglogaxis}[enlargelimits=false,xmin=6,xmax=3e9,ymin=3.5e-5,ymax=2,ytick={1e-5,1e-4,1e-3,1e-2,1e-1,1},height=0.4\textwidth, width=0.95\textwidth, grid=major, xlabel={number of samples $M$}, ylabel={rel\_{}sampl\_{}err\_A}, font=\footnotesize,
  legend style={at={(0.4,1.05)}, anchor=south,legend columns=4,legend cell align=left, font=\footnotesize, %
  },
  ]
\addplot[forget plot,dashed,mark options={solid},magenta,mark=pentagon,mark size=2.5] coordinates {
 (15,8.706e-01) (183,4.316e-01) (1643,2.125e-01) (12543,8.328e-02) (84845,2.895e-02) (574275,9.281e-03) (4068807,2.745e-03) (27910471,7.497e-04) (179044805,1.880e-04)
};
\addplot[magenta,mark=pentagon*,mark size=2.5] coordinates {
 (95,7.648e-01) (1415,4.166e-01) (9637,2.114e-01) (34613,9.100e-02) (144459,3.167e-02) (547947,1.065e-02) (2150373,3.219e-03) (7304585,8.897e-04) (29616599,2.227e-04) (82621003,5.304e-05)
};
\addlegendentry{$d$=7}
\addplot[forget plot,dashed,mark options={solid},red,mark=square,mark size=2.5] coordinates {
 (21,1.098e+00) (399,6.704e-01) (6753,4.195e-01) (78601,2.255e-01) (831125,1.075e-01) (7057695,4.612e-02) (69268743,1.803e-02)
};
\addplot[red,mark=square*,mark size=2.5] coordinates {
 (181,9.589e-01) (2269,6.649e-01) (22255,4.178e-01) (158907,2.383e-01) (908271,1.154e-01) (4123493,5.026e-02) (17137285,2.005e-02) (71759089,7.235e-03) (283874387,2.438e-03)
};
\addlegendentry{$d$=10}
\addplot[forget plot,dashed,mark options={solid},cyan,mark=triangle,mark size=2.5] coordinates {
 (31,1.198e+00) (1269,9.036e-01) (36015,7.202e-01) (683861,5.284e-01) (12971315,3.492e-01) (132560341,2.091e-01)
};
\addplot[cyan,mark=triangle*,mark size=2.5] coordinates {
 (197,1.067e+00) (4885,9.393e-01) (80655,7.359e-01) (800905,5.391e-01) (7323769,3.632e-01) (47110143,2.194e-01) (236566189,1.205e-01) (1224784507,6.081e-02)
};
\addlegendentry{$d$=15}
\addplot[forget plot,dashed,mark options={solid},darkgreen,mark=diamond,mark size=2.5] coordinates {
 (41,1.171e+00) (2641,9.959e-01) (121023,8.889e-01) (4249513,7.665e-01) (83869477,6.149e-01)
};
\addplot[darkgreen,mark=diamond*,mark size=2.5] coordinates {
 (267,1.071e+00) (10237,1.028e+00) (185353,9.288e-01) (3613339,7.822e-01) (28508819,6.268e-01) (219485937,4.693e-01) (1661989977,3.221e-01)
};
\addlegendentry{$d$=20}
\legend{};
  \end{loglogaxis}
\end{tikzpicture}
}
\vspace{1em}
\subfloat[$\Vert G_3^d - S^\Lambda_{I} G_3^d\vert L_2(\T^d)\Vert / \Vert G_3^d\vert L_2(\T^d)\Vert$
]{\label{fig:numerics:periodic:G3:rel_sampl_err_L2:M:hd}
\begin{tikzpicture}
  \begin{loglogaxis}[enlargelimits=false,xmin=6,xmax=3e9,ymin=1.5e-7,ymax=2e1,ytick={1e-6,1e-5,1e-4,1e-3,1e-2,1e-1,1},height=0.40\textwidth, width=0.95\textwidth, grid=major, xlabel={number of samples $M$}, ylabel={rel\_{}sampl\_{}err\_L2}, font=\footnotesize,
  legend style={at={(0.4,1.05)}, anchor=south,legend columns=4,legend cell align=left, font=\footnotesize, %
  },
  legend to name=samperrG3dMhd,
  ]
\addplot[forget plot,dashed,mark options={solid},magenta,mark=pentagon,mark size=2.5] coordinates {
 (15,5.044e-01) (183,8.648e-02) (1643,2.405e-02) (12543,5.702e-03) (84845,1.276e-03) (574275,2.684e-04) (4068807,5.291e-05) (27910471,9.777e-06) (179044805,1.689e-06)
};
\addplot[magenta,mark=pentagon*,mark size=2.5] coordinates {
 (95,3.262e-01) (1415,8.091e-02) (9637,2.357e-02) (34613,6.093e-03) (144459,1.315e-03) (547947,2.807e-04) (2150373,5.516e-05) (7304585,1.013e-05) (29616599,1.732e-06) (82621003,2.823e-07)
};
\addlegendentry{$d$=7}
\addplot[forget plot,dashed,mark options={solid},red,mark=square,mark size=2.5] coordinates {
 (21,1.165e+00) (399,1.486e-01) (6753,4.518e-02) (78601,1.331e-02) (831125,3.669e-03) (7057695,9.383e-04) (69268743,2.242e-04)
};
\addplot[red,mark=square*,mark size=2.5] coordinates {
 (181,6.360e-01) (2269,1.422e-01) (22255,4.475e-02) (158907,1.400e-02) (908271,3.811e-03) (4123493,9.705e-04) (17137285,2.323e-04) (71759089,5.174e-05) (283874387,1.094e-05)
};
\addlegendentry{$d$=10}
\addplot[forget plot,dashed,mark options={solid},cyan,mark=triangle,mark size=2.5] coordinates {
 (31,3.672e+00) (1269,2.868e-01) (36015,9.184e-02) (683861,3.383e-02) (12971315,1.160e-02) (132560341,3.691e-03)
};
\addplot[cyan,mark=triangle*,mark size=2.5] coordinates {
 (197,1.765e+00) (4885,3.840e-01) (80655,9.919e-02) (800905,3.480e-02) (7323769,1.202e-02) (47110143,3.806e-03) (236566189,1.126e-03) (1224784507,3.140e-04)
};
\addlegendentry{$d$=15}
\addplot[forget plot,dashed,mark options={solid},darkgreen,mark=diamond,mark size=2.5] coordinates {
 (41,1.028e+01) (2641,6.131e-01) (121023,1.487e-01) (4249513,6.317e-02) (83869477,2.521e-02)
};
\addplot[darkgreen,mark=diamond*,mark size=2.5] coordinates {
 (267,4.688e+00) (10237,9.718e-01) (185353,2.232e-01) (3613339,6.733e-02) (28508819,2.593e-02) (219485937,9.647e-03) (1661989977,3.318e-03)
};
\addlegendentry{$d$=20}
  \end{loglogaxis}
\end{tikzpicture}
}
\\[1em]
\pgfplotslegendfromname{samperrG3dMhd}
\caption{Relative $\mathcal{A}(\T^d)$
and $L_2(\T^d)$ sampling errors 
for $G_3^d$ with respect to the number of sampling nodes~$M$ for reconstructing single rank-1 lattices (dashed lines, unfilled markers) and reconstructing multiple rank-1 lattices (solid lines, filled markers),
when using the frequency index sets $I:=I_{N,\mathrm{even}}^{d,0}$.
}\label{fig:numerics:periodic:G3:rel_sampl_err:M:hd}
\end{figure}


\subsection{Kink test functions using dyadic hyperbolic cross frequency index sets}

In this section, we compare sampling along single and multiple rank-1 lattices with sparse grid sampling. For the latter, we use the \textsc{Matlab}$^{\text{\textregistered}}$ toolbox~\cite{DoKaKuPo_NHCFFT}.
Since this implementation uses a different type of hyperbolic cross index sets, we introduce additional notation and convert the theoretical results from Section~\ref{sec:mult_r1l_sampl}. The dyadic (non-symmetric) hyperbolic crosses are defined by
$H_n^d:=\bigcup_{\boldj\in\N_0^d,\; \|\boldj\|_1=n} Q_{\boldj}$, where $n\in\N_0$ denotes the refinement
and $Q_{\boldj}:=\bigtimes_{s=1}^d Q_{j_s}$, $Q_{j_s}:=\{1-2^{j_s-1},\ldots,2^{j_s-1}\}$.
When using these dyadic hyperbolic cross index sets~$H_n^d$ instead of $I_{N}^{d,0}$, one can show the same asymptotic upper bound as in \cite[Corollaries~2.40 and~2.44]{volkmerdiss} when using reconstructing single rank-1 lattices, see also \cite{ByKaUlVo16}.
Similarly, the results in Corollaries~\ref{cor:sampl_error_mr1l:Linf:Aabg} and~\ref{cor:approx_error:Hrtg:Habg} for reconstructing multiple rank-1 lattices
can be easily adapted.

Now, we approximate the scaled periodized (tensor product) kink function
\begin{equation} \label{equ:kink:25121}
g(\boldx):= \prod_{t=1}^d \frac{121 \sqrt{33}}{100} \max\left\{\frac{25}{121} - \left(x_t - \frac{1}{2}\right)^2\!, 0\right\}, \quad \boldx:=(x_1,\ldots,x_d)^\top\in\T^d,
\end{equation}
similar to \cite{HiMaOeUl15}.
We remark that we have $g\in\mathcal{A}^{1-\varepsilon}_\mathrm{mix}(\T^d)$ and $g\in\mathcal{H}^{3/2-\varepsilon}_\mathrm{mix}(\T^d)$, $\varepsilon>0$, as well as $\|g|\mathcal{A}(\T^d)\|\approx (1.84190)^d$ and $\|g|L_2(\T^d)\|=1$.

\begin{Example}
We use dyadic hyperbolic cross frequency index sets $I=H_n^d$ and approximate the kink function~$g$ based on samples along sparse grids and lattices.
For dimensions $d=3$ and $4$, we depict the obtained relative $\mathcal{A}(\T^d)$ and $L_2(\T^d)$ sampling errors in Figure~\ref{fig:kink25121:d3_d4}.
\begin{figure}[tb]
\centering
\subfloat[$d=3$]{
\begin{tikzpicture}
  \begin{loglogaxis}[font=\footnotesize,enlargelimits=false,xmin=1e0,xmax=8e7,ymin=1.5e-3,ymax=2e0,ytick={1e-3,1e-2,1e-1,1},height=0.4\textwidth, width=0.45\textwidth, grid=major, xlabel={number of samples $M$},
  ylabel={rel\_{}sampl\_{}err\_{}A},
  legend style={at={(0.5,1.05)}, anchor=south,legend columns=3},
  ]
\addplot[red,mark=x,mark size=2.5] coordinates {
 (4,1.000e+00) (13,1.000e+00) (38,1.001e+00) (104,7.968e-01) (272,4.882e-01) (688,2.660e-01) (1696,1.471e-01) (4096,9.052e-02) (9728,5.395e-02) (22784,3.056e-02) (52736,1.697e-02) (120832,9.228e-03) (274432,5.240e-03) (618496,3.003e-03)
};
\addlegendentry{SG}
\addplot[blue,dashed,mark options={solid},mark=diamond,mark size=2] coordinates {
 (4,9.374e-01) (17,6.803e-01) (85,4.233e-01) (325,2.725e-01) (1357,1.624e-01) (5065,9.767e-02) (17683,5.722e-02) (56663,3.334e-02) (217909,1.941e-02) (1045307,1.093e-02) (4118111,6.196e-03) (14619071,3.491e-03)
};
\addlegendentry{single}
\addplot[blue,mark=diamond*,mark size=2] coordinates { 
 (7,8.370e-01) (95,6.299e-01) (249,4.431e-01) (1119,3.142e-01) (2233,1.839e-01) (8463,1.109e-01) (20545,6.439e-02) (65777,3.707e-02) (214379,2.170e-02) (501895,1.221e-02) (1160759,7.016e-03) (2658889,3.902e-03)
};
 \addlegendentry{multiple}
\legend{};
  \end{loglogaxis}
\end{tikzpicture}
\hspace{2mm}
\hfill
\hspace{2mm}
\begin{tikzpicture}
  \begin{loglogaxis}[font=\footnotesize,enlargelimits=false,xmin=1e0,xmax=8e7,ymin=2e-5,ymax=2.5e0,ytick={1e-5,1e-4,1e-3,1e-2,1e-1,1},
  height=0.4\textwidth, width=0.45\textwidth, grid=major, xlabel={number of samples $M$},
  ylabel={rel\_{}sampl\_{}err\_{}L2},
  legend style={at={(0.5,1.05)}, anchor=south,legend columns=3},
  ]
\addplot[red,mark=x,mark size=2.5] coordinates {
 (4,1.000e+00) (13,1.000e+00) (38,8.875e-01) (104,5.362e-01) (272,2.224e-01) (688,7.011e-02) (1696,2.199e-02) (4096,8.742e-03) (9728,3.744e-03) (22784,1.374e-03) (52736,4.701e-04) (120832,1.652e-04) (274432,6.187e-05) %
};
\addlegendentry{SG}
\addplot[blue,dashed,mark options={solid},mark=diamond,mark size=2] coordinates {
 (4,8.475e-01) (17,4.463e-01) (85,2.343e-01) (325,1.123e-01) (1357,5.510e-02) (5065,2.102e-02) (17683,6.777e-03) (56663,2.734e-03) (217909,1.046e-03) (1045307,3.885e-04) (4118111,1.494e-04) (14619071,5.512e-05)
};
\addlegendentry{single}
\addplot[blue,mark=diamond*,mark size=2] coordinates { 
(7,6.620e-01) (95,3.778e-01) (249,2.416e-01) (1119,1.210e-01) (2233,5.794e-02) (8463,2.198e-02) (20545,7.064e-03) (65777,2.809e-03) (214379,1.072e-03) (501895,3.973e-04) (1160759,1.532e-04) (2658889,5.631e-05)
};
 \addlegendentry{multiple}
\addplot[black,domain=4e4:7e5,samples=100,dashdotted,very thick] {1.1*x^(-1.5)*ln(x)^((3-1)*(1.5+1/2))};
\addlegendentry{theor.}
\addplot[black,domain=1e5:4e7,samples=100,dotted,very thick] {0.4*x^(-1.5/2)*ln(x)^((3-2)*1.5/2+(3-1)/2)};
 \addlegendentry{theor.}
\addplot[black,domain=1e5:7e6,samples=100,thick] {0.005*x^(-1)*ln(x)^(3*1+1)};
\addlegendentry{theor.}
\legend{};
  \end{loglogaxis}
\end{tikzpicture}
}
\\[1.5em]
\subfloat[$d=4$]{
\begin{tikzpicture}
  \begin{loglogaxis}[font=\footnotesize,enlargelimits=false,xmin=1e0,xmax=8e7,ymin=7e-3,ymax=2e0,ytick={1e-2,1e-1,1},height=0.4\textwidth, width=0.45\textwidth, grid=major, xlabel={number of samples $M$},
  ylabel={rel\_{}sampl\_{}err\_{}A},
  legend style={at={(0.5,1.05)}, anchor=south,legend columns=3},
  ]
\addplot[red,mark=x,mark size=2.5] coordinates {
 (5,1.000e+00) (19,1.000e+00) (63,1.000e+00) (192,1.052e+00) (552,9.897e-01) (1520,7.520e-01) (4048,4.810e-01) (10496,2.868e-01) (26624,1.785e-01) (66304,1.127e-01) (162560,6.900e-02) (393216,4.049e-02) (940032,2.339e-02) (2224128,1.351e-02)
};
\addlegendentry{SG}
\addplot[blue,dashed,mark options={solid},mark=diamond,mark size=2] coordinates {
 (5,9.156e-01) (36,7.902e-01) (191,6.203e-01) (724,4.396e-01) (3420,3.018e-01) (16936,1.973e-01) (59719,1.268e-01) (288420,7.933e-02) (1199400,4.893e-02) (4880078,2.970e-02) (28876381,1.787e-02)
};
\addlegendentry{single}
 \addplot[blue,mark=diamond*,mark size=2] coordinates { 
 (23,9.203e-01) (37,8.198e-01) (393,6.223e-01) (3243,4.758e-01) (7879,3.340e-01) (30809,2.194e-01) (81459,1.398e-01) (210293,8.986e-02) (533057,5.517e-02) (1591757,3.326e-02) (3902341,2.015e-02) (9438331,1.201e-02)
 };
 \addlegendentry{multiple}
\legend{};
  \end{loglogaxis}
\end{tikzpicture}
\hspace{2mm}
\hfill
\hspace{2mm}
\begin{tikzpicture}
  \begin{loglogaxis}[font=\footnotesize,enlargelimits=false,xmin=1e0,xmax=8e7,ymin=7e-5,ymax=2.5e0,ytick={1e-4,1e-3,1e-2,1e-1,1},height=0.4\textwidth, width=0.45\textwidth, grid=major, xlabel={number of samples $M$},
  ylabel={rel\_{}sampl\_{}err\_{}L2},
  legend style={at={(0.5,1.05)}, anchor=south,legend columns=3,legend cell align=left,font=\footnotesize},
  legend to name=samperrkink,
  ]
\addplot[red,mark=x,mark size=2.5] coordinates {
 (5,1.000e+00) (19,1.000e+00) (63,1.000e+00) (192,9.972e-01) (552,7.362e-01) (1520,3.916e-01) (4048,1.585e-01) (10496,5.407e-02) (26624,1.873e-02) (66304,7.645e-03) (162560,3.187e-03) (393216,1.187e-03) (940032,4.179e-04) (2224128,1.510e-04)
};
\addlegendentry{SG}
\addplot[blue,dashed,mark options={solid},mark=diamond,mark size=2] coordinates {
 (5,7.995e-01) (36,5.115e-01) (191,3.331e-01) (724,1.857e-01) (3420,1.063e-01) (16936,5.106e-02) (59719,2.362e-02) (288420,9.619e-03) (1199400,3.493e-03) (4880078,1.382e-03) (28876381,5.468e-04)
};
\addlegendentry{single}
 \addplot[blue,mark=diamond*,mark size=2] coordinates { 
 (23,7.634e-01) (37,5.531e-01) (393,3.341e-01) (3243,1.923e-01) (7879,1.104e-01) (30809,5.226e-02) (81459,2.413e-02) (210293,9.867e-03) (533057,3.572e-03) (1591757,1.407e-03) (3902341,5.570e-04) (9438331,2.130e-04)
 };
 \addlegendentry{multiple}
\addplot[black,domain=3e5:5e6,samples=100,dashdotted,very thick] {0.14*x^(-1.5)*ln(x)^((4-1)*(1.5+1/2))};
\addlegendentry{theor. SG}
\addplot[black,domain=3e5:4e7,samples=100,dotted,very thick] {0.1*x^(-1.5/2)*ln(x)^((4-2)*1.5/2+(4-1)/2)};
\addlegendentry{theor. single}
\addplot[black,domain=1e5:2e7,samples=100,thick] {0.005*x^(-1)*ln(x)^(4*1+1)};
\addlegendentry{theor. mult.}
  \end{loglogaxis}
\end{tikzpicture}
}
\\[1em]
\pgfplotslegendfromname{samperrkink}
\caption{Relative $\mathcal{A}(\T^d)$ and $L_2(\T^d)$ sampling errors for the approximation of the kink function~$g$ from~\eqref{equ:kink:25121} when sampling along sparse grids (SG) as well as single and multiple rank-1 lattices.}\label{fig:kink25121:d3_d4}
\end{figure}
The obtained errors decay for increasing refinements~$n$ of the hyperbolic crosses~$H_n^d$ and consequently larger number of sampling points~$M$.
For the considered test function~$g$ in dimensions $d=3$ and $4$, the sparse grid sampling performs slightly better than the lattice sampling.
Additionally, we plot the theoretical asymptotic upper bounds for the relative $L_2(\T^d)$ sampling errors
for sparse grids, which are $\mathcal{O}\big(M^{-\beta} (\log M)^{(d-1)(\beta+1/2)}\big)$ due to \cite[Theorem~6.10]{ByDuSiUl14}, as dash-dotted graphs,
for single rank-1 lattices, which are $\mathcal{O}\big(M^{-\beta/2} \, (\log M)^{\frac{d-2}{2}\beta+\frac{d-1}{2}}\big)$ due to \cite[Corollary~2.44]{volkmerdiss}, as dotted lines as well as
for multiple rank-1 lattices, which are $\mathcal{O}\big(M^{-(\beta-\lambda)} \, (\log M)^{(\beta-\lambda)d+1}\big)$ due to Corollary~\ref{cor:approx_error:Hrtg:Habg}, as solid lines without markers,
where $\beta=3/2-\epsilon$, $\epsilon>0$, and $\lambda>1/2$ in the exponent.
We observe that the errors for sparse grid and single rank-1 lattice sampling approximately decay as the upper bounds suggest.
For multiple rank-1 lattice sampling, the obtained errors seem to decay faster by an additive factor of about 1/2 in the main term, i.e., we may not observe the term $\lambda>1/2$.
\newline
Additionally, we consider the relative $L_2(\T^d)$ sampling errors in dimensions $d=6,\ldots,9$ and visualize the error values in Figure~\ref{fig:kink25121:d6_d9}. Here, we observe for increasing dimension~$d$ that the sparse grids yield a worse pre-asymptotic behavior compared to the single and multiple rank-1 lattices. For instance for $d=9$, the relative $L_2(\T^d)$ sampling errors are still greater than 1 when using approximately 10 million samples, whereas the corresponding errors are smaller by one order of magnitude when using single or multiple rank-1 lattices.
\end{Example}

\begin{figure}[tb]
\centering
\subfloat[$d=6$]{
\begin{tikzpicture}
  \begin{loglogaxis}[font=\footnotesize,enlargelimits=false,xmin=3e0,xmax=3e8,ymin=1e-3,ymax=3e0,ytick={1e-4,1e-3,1e-2,1e-1,1},height=0.4\textwidth, width=0.45\textwidth, grid=major, xlabel={number of samples $M$},
   ylabel={rel\_{}sampl\_{}err\_{}L2},
  legend style={at={(0.5,1.05)}, anchor=south,legend columns=3},
  ]
\addplot[red,mark=x,mark size=2.5] coordinates {
 (7,1.000e+00) (34,1.000e+00) (138,1.000e+00) (501,1.000e+00) (1683,1.000e+00) (5336,1.160e+00) (16172,1.132e+00) (47264,8.347e-01) (134048,4.802e-01) (370688,2.253e-01) (1003136,9.110e-02) (2664192,3.447e-02) (6960384,1.355e-02) (17922048,5.629e-03) (45552640,2.280e-03)
};
\addlegendentry{SG}
\addplot[blue,dashed,mark options={solid},mark=diamond,mark size=2] coordinates {
 (7,8.384e-01) (83,7.360e-01) (699,4.936e-01) (3988,3.347e-01) (19719,2.244e-01) (127882,1.365e-01) (552353,8.035e-02) (2954249,4.334e-02) (13539013,2.230e-02)
};
\addlegendentry{single}
 \addplot[blue,mark=diamond*,mark size=2] coordinates { 
 (29,8.329e-01) (287,6.583e-01) (1131,5.004e-01) (6121,3.502e-01) (27219,2.329e-01) (118017,1.403e-01) (323789,8.263e-02) (1040157,4.433e-02) (3218403,2.289e-02) (10380567,1.081e-02) (30095325,4.874e-03) (85256153,2.079e-03)
 };
 \addlegendentry{multiple}
\addplot[black,domain=2e6:1e8,samples=100,dashdotted,very thick] {0.00013*x^(-1.5)*ln(x)^((6-1)*(1.5+1/2))};
\addlegendentry{theor.}
\addplot[black,domain=1e6:1e8,samples=100,dotted,very thick] {0.003*x^(-1.5/2)*ln(x)^((6-2)*1.5/2+(6-1)/2)};
\addlegendentry{theor.}
\addplot[black,domain=1e7:2e8,samples=100,thick] {0.0006*x^(-1)*ln(x)^(6*1+1)};
\addlegendentry{theor.}
\legend{};
  \end{loglogaxis}
\end{tikzpicture}
}
\hfill
\subfloat[$d=7$]{
\begin{tikzpicture}
  \begin{loglogaxis}[font=\footnotesize,enlargelimits=false,xmin=3e0,xmax=3e8,ymin=8e-3,ymax=2e0,ytick={1e-4,1e-3,1e-2,1e-1,1},height=0.4\textwidth, width=0.45\textwidth, grid=major, xlabel={number of samples $M$},
  ylabel={rel\_{}sampl\_{}err\_{}L2},
  legend style={at={(0.5,1.05)}, anchor=south,legend columns=3,font=\footnotesize},
  legend to name=samperrkinkhd,
  ]
\addplot[red,mark=x,mark size=2.5] coordinates {
 (8,1.0) (43,1.0) (190,1.0) (743,1.0) (2668,1.0) (8989,1.0) (28814,1.223) (88720,1.329) (264224,1.108) (765088,7.252e-01) (2162624,3.881e-01) (5986304,1.769e-01) (16268800,7.226e-02) (43499264,2.846e-02) (114629120,1.157e-02)
};
\addlegendentry{SG}
\addplot[blue,dashed,mark options={solid},mark=diamond,mark size=2] coordinates {
 (8,9.284e-01) (127,7.358e-01) (979,5.739e-01) (7717,4.065e-01) (46879,2.864e-01) (252890,1.871e-01) (1202231,1.180e-01) (7895149,6.927e-02) (49831729,3.897e-02)
};
\addlegendentry{single}
\addplot[blue,mark=diamond*,mark size=2] coordinates {
 (57,9.584e-01) (387,7.283e-01) (2353,5.713e-01) (12069,4.199e-01) (53905,2.940e-01) (180177,1.917e-01) (519013,1.208e-01) (1952717,7.086e-02) (6342695,3.984e-02) (22954985,2.091e-02) (73530869,1.048e-02)
 };
 \addlegendentry{multiple}
\addplot[black,domain=2e6:1e9,samples=100,dashdotted,very thick] {0.00002*x^(-1.5)*ln(x)^((7-1)*(1.5+1/2))};
\addlegendentry{theor. SG}
\addplot[black,domain=1e5:2e8,samples=100,dotted,very thick] {0.00014*x^(-1.5/2)*ln(x)^((7-2)*1.5/2+(7-1)/2)};
\addlegendentry{theor. single}
\addplot[black,domain=5e6:7e7,samples=100,thick] {0.000095*x^(-1)*ln(x)^(7*1+1)};
\addlegendentry{theor. mult.}
  \end{loglogaxis}
\end{tikzpicture}
}
\\
\subfloat[$d=8$]{
\begin{tikzpicture}
  \begin{loglogaxis}[font=\footnotesize,enlargelimits=false,xmin=3e0,xmax=3e8,ymin=1e-2,ymax=2e0,ytick={1e-4,1e-3,1e-2,1e-1,1},height=0.4\textwidth, width=0.45\textwidth, grid=major, xlabel={number of samples $M$},
  ylabel={rel\_{}sampl\_{}err\_{}L2},
  legend style={at={(0.5,1.05)}, anchor=south,legend columns=3},
  ]
\addplot[red,mark=x,mark size=2.5] coordinates {
 (9,1.0) (53,1.0) (253,1.0) (1059,1.0) (4043,1.0) (14407,1.0) (48639,1.0) (157184,1.278e+00) (489872,1.528e+00) (1480608,1.416e+00) (4358752,1.035e+00) (12541184,6.205e-01)
};
\addlegendentry{SG}
\addplot[blue,dashed,mark options={solid},mark=diamond,mark size=2] coordinates {
 (9,9.303e-01) (160,8.303e-01) (1603,6.317e-01) (13027,4.840e-01) (84572,3.482e-01) (657990,2.406e-01) (2988071,1.602e-01) (20676791,1.004e-01) (160169101,6.043e-02) 
};
\addlegendentry{single}
 \addplot[blue,mark=diamond*,mark size=2] coordinates {
 (57,9.800e-01) (583,8.455e-01) (2637,6.464e-01) (14989,4.881e-01) (73021,3.574e-01) (288587,2.464e-01) (1168279,1.632e-01) (4402667,1.022e-01) (13718045,6.158e-02) (44420491,3.497e-02) (139482127,1.905e-02)
};
 \addlegendentry{multiple}
\legend{};
  \end{loglogaxis}
\end{tikzpicture}
}
\hfill
\subfloat[$d=9$]{
\begin{tikzpicture}
  \begin{loglogaxis}[font=\footnotesize,enlargelimits=false,xmin=3e0,xmax=3e8,ymin=1e-2,ymax=3e0,ytick={1e-4,1e-3,1e-2,1e-1,1},height=0.4\textwidth, width=0.45\textwidth, grid=major, xlabel={number of samples $M$},
  ylabel={rel\_{}sampl\_{}err\_{}L2},
  legend style={at={(0.5,1.05)}, anchor=south,legend columns=3,legend cell align=left},
  ]
\addplot[red,mark=x,mark size=2.5] coordinates {
 (10,1.000e+00) (64,1.000e+00) (328,1.000e+00) (1462,1.000e+00) (5908,1.000e+00) (22180,1.000e+00) (78592,1.000e+00) (265729,1.000e+00) (864146,1.326e+00) (2719028,1.730e+00) (8316200,1.761e+00) (24814832,1.418e+00)
};
\addlegendentry{SG}
\addplot[blue,dashed,mark options={solid},mark=diamond,mark size=2] coordinates {
 (10,9.724e-01) (212,8.977e-01) (2878,6.847e-01) (21058,5.376e-01) (199896,4.085e-01) (1460690,2.953e-01) (6637979,2.057e-01) (52963709,1.359e-01) %
};
\addlegendentry{single}
 \addplot[blue,mark=diamond*,mark size=2] coordinates {
 (69,9.296e-01) (679,8.176e-01) (4039,6.907e-01) (23671,5.527e-01) (106657,4.199e-01) (444305,3.026e-01) (2201313,2.092e-01) (7441829,1.384e-01) (25926193,8.792e-02) (87010369,5.299e-02)
 };
 \addlegendentry{multiple}
\legend{};
  \end{loglogaxis}
\end{tikzpicture}
}
\\[1em]
\pgfplotslegendfromname{samperrkinkhd}
\caption{Relative $L_2(\T^d)$ sampling errors for the approximation of the kink function~$g$ from~\eqref{equ:kink:25121} when sampling along sparse grids (SG) as well as single and multiple rank-1 lattices.}\label{fig:kink25121:d6_d9}
\end{figure}

\section*{Acknowledgements}
We thank the referees for the valuable suggestions.
LK gratefully acknowledges the funding by the Deutsche Forschungsgemeinschaft 
(DFG, German Research Foundation, project number 380648269).
TV gratefully acknowledges the funding by the European Union and the Free State of Saxony (EFRE/ESF).

\begin{small}
\bibliographystyle{abbrv}

\end{small}

\end{document}